\documentclass[11pt]{article}%
\pdfoutput=1 %
\usepackage{amssymb}
\usepackage[margin=1in]{geometry}

\usepackage{lineno}
\usepackage{graphicx}
\usepackage{moreverb,url}
\usepackage[ruled,vlined]{algorithm2e}
\usepackage[colorlinks,bookmarksopen,bookmarksnumbered,citecolor=red,urlcolor=red]{hyperref}

\usepackage{xcolor,tabularx,nccmath}
\usepackage{amsmath,amsthm,enumitem}
\usepackage{booktabs,caption,subcaption,mathtools,multirow}
\newtheorem{definition}{Definition}
\newcommand{\bigM}{\mathbb{M}}
\DeclareMathAlphabet{\mathdutchcal}{U}{dutchcal}{m}{n}
\SetMathAlphabet{\mathdutchcal}{bold}{U}{dutchcal}{b}{n}
\DeclareMathAlphabet{\mathdutchbcal}{U}{dutchcal}{b}{n}
\newcommand{\sA}{\mathdutchcal{A}}
\newcommand{\sB}{\mathdutchcal{B}}
\newcommand{\sN}{\mathdutchcal{N}}
\newcommand{\sK}{\mathdutchcal{K}}

\newcommand{\sT}{\mathdutchcal{T}}
\newcommand{\sF}{\mathdutchcal{F}}

\newcommand{\sL}{\mathdutchcal{L}}

\newcommand{\sR}{\mathdutchcal{R}}
\newcommand{\sC}{\mathdutchcal{C}}

\newcommand{\sZ}{\mathdutchcal{Z}}

\usepackage[nameinlink,capitalise]{cleveref}
\hypersetup{colorlinks={true},allcolors={blue}}

\crefformat{section}{\S#2#1#3}
\crefformat{subsection}{\S#2#1#3}
\crefformat{equation}{#2Equation~\textup{(#1)}#3}
\crefformat{cons}{#2Constraint~\textup{(#1)}#3}
\labelcrefformat{cons}{\textup{#2(#1)#3}}
\crefformat{consset}{#2Constraints~\textup{(#1)}#3} %

\labelcrefformat{consset}{\textup{#2(#1)#3}}
\crefformat{func}{#2Function~\textup{(#1)}#3}
\labelcrefformat{func}{#2Function~\textup{(#1)}#3}
\crefformat{algo}{#2Algorithm~\textup{#1}#3}
\crefformat{objfunc}{#2Objective function~\textup{(#1)}#3}
\labelcrefformat{objfunc}{\textup{#2(#1)#3}}
\crefformat{ch}{#2Chapter~\textup{#1}#3}
\crefformat{sec}{#2Section~\textup{#1}#3}
\crefformat{tab}{#2Table~\textup{#1}#3}
\crefformat{fig}{#2Figure~\textup{#1}#3}
\crefdefaultlabelformat{#2#1#3} %
\Crefname{figure}{Figure}{Figures}
\crefformat{eq}{#2equation~\textup{(#1)}#3}
\labelcrefformat{eq}{\textup{(#2#1#3)}}
\crefformat{theo}{#2Theorem~\textup{#1}#3}
\crefformat{lem}{#2Lemma~\textup{#1}#3}
\crefformat{def}{#2Definition~\textup{#1}#3}
\crefformat{app}{#2\textup{\!#1}#3}
\crefformat{scenario}{#2Scenario~\textup{#1}#3}
 \AtBeginDocument{%
    \crefname{figure}{Figure}{figures}%
}
\let\origref\cref
\def\cref#1{\origref{#1}}
\renewenvironment{proof}[1][\proofname]{{\bfseries #1.}}{\qed}

\expandafter\def\expandafter\UrlBreaks\expandafter{\UrlBreaks%
  \do\a\do\b\do\c\do\d\do\e\do\f\do\g\do\h\do\i\do\j%
  \do\k\do\l\do\m\do\n\do\o\do\p\do\q\do\r\do\s\do\t%
  \do\u\do\v\do\w\do\x\do\y\do\z\do\A\do\B\do\C\do\D%
  \do\E\do\F\do\G\do\H\do\I\do\J\do\K\do\L\do\M\do\N%
  \do\O\do\P\do\Q\do\R\do\S\do\T\do\U\do\V\do\W\do\X%
  \do\Y\do\Z}

\usepackage{stackengine}
\newcommand\barbelow[1]{\stackunder[1.2pt]{$#1$}{\rule{.8ex}{.075ex}}}
\usepackage{algpseudocode}

\DeclarePairedDelimiter\ceil{\lceil}{\rceil}

\newcommand{\norm}[1]{\lvert #1 \rvert}
\usepackage[normalem]{ulem} %
 
\newcommand{\replace}[2]{{\color{black!50!black}\color{black}{\color{black}#2\color{black}}}} %

\newtheorem{theorem}{Theorem}%
\newtheorem{lemma}[theorem]{Lemma}

\crefname{defi}{definition}{definitions}
\Crefname{defi}{Definition}{Definitions}
\crefname{lemma}{lemma}{lemmas}
\Crefname{lemma}{Lemma}{Lemmas}
\crefname{assumption}{assumption}{assumptions}
\Crefname{assumption}{Assumption}{Assumptions}

\usepackage{multicol}

\usepackage{authblk}
\providecommand{\keywords}[1]{\noindent\textbf{Keywords:} #1}

\usepackage{natbib}

\begin{document}

\title{Electric Vehicle Supply Equipment Location and Capacity Allocation for Fixed-Route Networks}

\author[1]{Amir Davatgari}

\affil[1]{University of Illinois Chicago, 1200 W. Harrison St., Chicago, 60607, IL, USA}

\author[2,3]{Taner Cokyasar}
\author[2,4]{Anirudh Subramanyam}
\author[2]{Jeffrey Larson}
\author[1]{Abolfazl (Kouros) Mohammadian}

\affil[2]{Argonne National Laboratory, 9700 S. Cass Avenue, Lemont, 60439, IL, USA}
\affil[3]{TrOpt R\&D, Balcali mah., Saricam, 01330, Adana, Turkey}
\affil[4]{The Pennsylvania State University, 201 Old Main, University Park, 16802, PA, USA}

\maketitle

\begin{abstract}
Electric vehicle (EV) supply equipment location and allocation (EVSELCA) problems for freight vehicles are becoming more important because of the trending electrification shift. Some previous works address EV charger location and vehicle routing problems simultaneously by generating vehicle routes from scratch. Although such routes can be efficient, introducing new routes may violate practical constraints, such as drive schedules, and satisfying electrification requirements can require dramatically altering existing routes. To address the challenges in the prevailing adoption scheme, we approach the problem from a fixed-route perspective. We develop a mixed-integer linear program, a clustering approach, and a metaheuristic solution method using a genetic algorithm (GA) to solve the EVSELCA problem. The clustering approach simplifies the problem by grouping customers into clusters, while the GA generates solutions that are shown to be nearly optimal for small problem cases. A case study examines how charger costs, energy costs, the value of time (VOT), and battery capacity impact the cost of the EVSELCA. Charger equipment costs were found to be the most significant component in the objective function, leading to a substantial reduction in cost when decreased. VOT costs exhibited a significant decrease with rising energy costs. An increase in VOT resulted in a notable rise in the number of fast chargers. Longer EV ranges decrease total costs up to a certain point, beyond which the decrease in total costs is negligible.
\end{abstract}

\keywords{
EV charging facility location and capacity allocation, freight electrification, truck electrification, fixed-route electrification, optimization}

\section{Introduction}
The widespread use of fossil fuels to meet energy requirements produces greenhouse gas (GHG) emissions that negatively impact the climate and environment \citep{Metz2007}. According to the \cite{USEnvironmentalProtectionAgency2020}, the transportation sector is the most significant contributor to GHG emissions in the United States, accounting for 27\% of the total emissions. Medium-duty (MD) and heavy-duty (HD) trucks are responsible for 26\% of the total emissions in this sector. With the increase in e-commerce, the number of trucks and their total miles traveled are expected to grow,  resulting in higher emissions from freight transportation \citep{HovlandConsultingLLC2020}. 

Truck electrification is a promising solution that can help mitigate the negative impact of their GHG emissions \citep{TALEBIAN2018109}. By replacing conventional engines with electric motors, the transportation sector can significantly reduce GHG emissions. While conventional trucks can cover long distances without refueling, a critical factor for electric trucks is their relatively shorter battery range \citep{HovlandConsultingLLC2020}.  The range of a typical MD/HD electric truck is around 130 miles \citep{Battery_capacity}, which limits their use in long-haul trucking. Additionally, the time required for recharging electric trucks can be another significant challenge: Depending on the charging method, it can take 20 minutes to 8 hours to fully recharge their batteries \citep{Bennett2021}. This can lead to significant downtimes for trucking companies, which affect their productivity and profitability. Therefore, efficiently solving electric vehicle supply equipment location and capacity allocation (EVSELCA) problems is crucial to making electric trucks a viable option for commercial transportation. The EVSELCA problem is to find optimal locations, numbers, and types of electric vehicle (EV) supply equipment (i.e., \textit{EV chargers} or \textit{chargers} for short) to minimize strategic investment costs while satisfying operational constraints.

We present a new approach to solve the EVSELCA problem by developing a mixed-integer linear programming (MILP) model that optimizes the locations and numbers of various types of chargers. The objective of the MILP model is to minimize strategic investment costs while satisfying operational constraints. To achieve this, our MILP model takes into account fixed-facility costs, charger costs, recharging energy costs, and value of time (VOT) costs. VOT costs account for the time spent traveling to a recharging station, waiting and recharging, and returning to service. By finding the optimal balance between these cost components, our model can aid in long-term EVSELCA planning. Because of the problem's complexity, we use a clustering approach that groups customers into clusters and allows recharging only after servicing these clusters. Moreover, we propose a metaheuristic solution method based on a genetic algorithm (GA) to generate near-optimal solutions within a reasonable time, enabling the model to be applied to large-scale instances. 

Our study contributes to the EVSELCA literature through four crucial aspects. First, when designing EVSELCA for the freight transportation industry, it is important to consider fixed routes \citep{GHAMAMI2016389,WANG2016242}. The EVSELCA problem has often been modeled as an EV location routing problem (EVLRP), which solves both the strategic charging facility location (and allocation in some cases) problem and the routing problem \citep{yang2015battery,hof2017solving,schiffer2017electric,schiffer2018designing,schiffer2018adaptive,RAEESI202282}. Although creating new routes for EVs can lead to better solutions \citep{Shojaei2022}, the convention in this industry is to use a fixed-route approach, usually electrifying existing routes that are shorter than an EV range. Second, the location and allocation decisions should be made jointly, with the numbers and types of chargers for each location serving as decision variables \citep{Davatgari2021,GHAMAMI2020102802,arias}. The type of charger impacts recharging time and infrastructure costs, while the number of chargers affects installation costs and waiting times for recharging. Thus, strategic planning requires a balance between waiting costs, recharging time, and infrastructure costs through the selection of appropriate types and numbers of chargers. Third, it is crucial to take into account the dynamic nature of charging demand over time since it plays a critical role, as highlighted by \cite{GHAMAMI2016389} and \cite{WANG2016242}. If multiple recharging events happen simultaneously, the design would require an excessive number of chargers. By considering dynamic charging demand for recharging, however, it becomes possible to schedule these events in a way that uses fewer chargers, resulting in higher utilization rates. Therefore, the spatio-temporal aspect of the problem is conserved. Fourth, partial recharging is crucial and should be considered in the EVSELCA problem \citep{LI2016128}. In some cases, electric trucks may be partially recharged because of operational time limitations, and a model without this consideration may produce impractical solutions. Overall, the main contribution of our study is its comprehensive consideration of these four key aspects of the EVSELCA problem. While other studies in the literature have addressed these aspects, they modeled each aspect either individually or as a combinatorial subset that lacked one of the other aspects.

In \cref{lit_rev} we first review studies that have used various approaches to tackle the EVSELCA problem. We provide a more comprehensive definition of the EVSELCA problem and clearly demonstrate the MILP in \cref{methodology}. Metaheuristic solution procedures are developed in \cref{heuristic}. The setup and results of numerical experiments are presented in \cref{case_studies} using data from POLARIS, the Planning and Operations Language for Agent-based Regional Integrated Simulation developed at Argonne National Laboratory \citep{auld2016polaris}. In \cref{conclusion} we conclude the study and discuss potential future research directions.

\section{Literature review} \label[sec]{lit_rev}
Many studies explore the strategic planning of EV charger placement \citep{GHAMAMI2020102802,ZHU201611,GHAMAMI2016389,LI2016128,WANG2016242,Davatgari2021,whitehead,LIU2018748,worley2012simultaneous,Speth_2022,arias,Kavianipour2021,Singh2022,Kavianipour}. Most of these studies focus on light-duty (LD) vehicles rather than MD and HD vehicles (hereafter called trucks) \citep{GHAMAMI2020102802,ZHU201611,LIU2018748,GHAMAMI2016389,WANG2016242,Kavianipour2021,Singh2022,Kavianipour}. For example, \cite{Davatgari2021} develops a mixed-integer nonlinear programming (MINLP) model to solve the EVSELCA problem considering routing for LD vehicles. Our study, in contrast, aims to solve the EVSELCA problem for trucks with fixed routes.  

The planning of EV chargers for trucks is different from that of LD vehicles because trucks often require fast recharging due to the high value of time in the business world and because their larger batteries necessitate longer recharging times. Although fast-charging equipment can reduce recharging time, it is also more costly. Additionally, unlike LD vehicles, trucks typically have predetermined routes and operate with time limits enforced by law and regulations \citep{operation_time}. Given the differences, the EVSELCA problem for trucks is an area that has not been extensively explored in the literature, and this study aims to address this gap. The studies conducted on this problem generally can be classified into two categories based on their methodology: (1) a coverage-oriented approach \citep{whitehead,Speth_2022} and (2) a demand-oriented approach \citep{worley2012simultaneous,arias,LIU2018748}; our study falls into the second category.

Coverage-oriented approaches aim to maximize the geographical coverage of a recharging service. This approach often assumes that each charger location can meet the recharging demand of a circular area. The objective is then to maximize the coverage of a region while minimizing the number of circles, which represents the number of charger locations. Such approaches often do not consider demand intensity, the potential impact of queuing during recharging events, and other operational constraints such as the inability to recharge two vehicles simultaneously using one charger. For instance, \citet{Speth_2022} model a coverage-oriented approach as a linear 
programming (LP) model to determine charger locations in order to minimize the number of chargers while maximizing the geographic coverage. After locating a charger, the model uses queuing to estimate the number of chargers in each EV charging facility. In contrast, our research focuses on locating EV chargers and optimizing the number of chargers using a demand-oriented approach that considers the impact of queuing during charging events.  

Demand-oriented approaches aim to minimize strategic (e.g., infrastructure and charger) costs  and operational (e.g., time and energy) costs subject to demand-conservation constraints that guarantee a certain level of service based on a deployment decision. This approach is often modeled as an EVLRP in the literature, which involves determining the allocation and location of EV chargers and solving the routing problem \citep{yang2015battery,hof2017solving,schiffer2017electric,schiffer2018designing,schiffer2018adaptive}. For example, \citet{worley2012simultaneous} develop an MILP to determine charger locations that minimize travel time, recharging costs, and construction costs. While our study shares a similar objective, our approach differs in the vehicle routing aspect because we focus on fixed routes. Additionally, our study takes into account several important factors that were not considered in the aforementioned research, such as allocation, dynamic charging demand, and partial recharging. By taking into account these additional variables, we aim to provide a more accurate and comprehensive solution. Another study by \citet{LIU2018748} uses a bi-level approach to model the EVLRP. In this approach, the upper level focuses on minimizing GHG emissions by determining the optimal location of chargers, while the lower level solves a mixed-traffic assignment of LD vehicles and trucks.  Our study differs from this approach as well. Specifically, we use a fixed-route approach to minimize the costs associated with charger placement, location, energy consumption, the value of time for recharging, and detouring. Additionally, like our previous example, we take into account allocation, dynamic charging demand, and partial recharging. \cref{literature} summarizes relevant studies and compares our study with their objectives, model types, and other features.

\begin{table}[!htb]
\scriptsize
  \setlength\extrarowheight{4pt}
  \caption{\replace{}{Summary of the existing relevant literature}}\label[tab]{literature}
  \resizebox{\textwidth}{!}{\begin{tabular}{p{1.3cm}p{0.8cm}>{\raggedright}p{4cm}p{.8cm}p{.8cm}p{.2cm}p{.3cm}p{.3cm}p{.3cm}p{.3cm}p{.4cm}p{.4cm}p{.4cm}}
  \toprule
  Study & Vehicle & Objective & Model & Method &(i)&(ii)&(iii)&(iv)&(v)&(vi)&(vii)&(viii)\\
  \midrule
\citet{GHAMAMI2020102802} & LD & Minimize infrastructure cost and users' detour, waiting, and charging delay & MINLP & SA & \checkmark & \checkmark & \checkmark & - & - & - & \checkmark & N/A \\
\citet{ZHU201611} & LD & Minimize construction costs and station access cost & MILP & GA & \checkmark & \checkmark & \checkmark & - & - & - & - & N/A \\
\citet{GHAMAMI2016389} & LD & Minimize infrastructure and battery costs, and recharging and queueing time & MINLP & SA & - & \checkmark & \checkmark & - & \checkmark & - & \checkmark & - \\
\citet{LI2016128} & LD & Minimize cost of construction and relocation of existing chargers & MILP & GA & \checkmark & \checkmark & - & - & - & \checkmark & - & N/A \\
\citet{WANG2016242} & LD & Minimize operational and construction costs & MILP & CS & - & \checkmark & - & - & \checkmark & - & - & \checkmark \\
\citet{Davatgari2021} & LD & Minimize total system travel time and construction cost of EV charging infrastructure & MINLP & GA & \checkmark & \checkmark & \checkmark & \checkmark & - & - & - & N/A \\
\citet{Kavianipour2021} & LD & Minimize the total system cost including charging station and charger installation costs, and charging, queuing, and detouring delays & MINLP & SA & \checkmark & \checkmark & \checkmark & - & \checkmark & \checkmark & \checkmark & N/A \\
\citet{whitehead} & MD & Maximize coverage & MILP & NS & \checkmark & \checkmark & - & - & - & - & - & N/A \\
\citet{LIU2018748} & MD & Minimize emissions & MPCC & H & \checkmark & \checkmark & - & - & - & - & - & N/A \\
\citet{worley2012simultaneous} & HD & Minimize transportation, recharging, and charging station placement costs & MILP & NS & \checkmark & \checkmark & - & - & - & - & - & N/A \\
\citet{Speth_2022} & HD & Maximize coverage & LP & H & \checkmark & \checkmark & - & - & - & - & - & N/A \\
\citet{arias} & HD & Minimize energy consumption, charger installation, and routing costs & MILP & CS & \checkmark & \checkmark & - & - & - & - & - & N/A \\
This study & MD/HD & Minimize charger, location, energy, value of time for recharging and detouring costs & MILP & GA & - & \checkmark & \checkmark & \checkmark & \checkmark & \checkmark & \checkmark & \checkmark \\ \bottomrule
\end{tabular}%
}
\replace{}{(i) Route choice,
(ii) Location planning,
(iii) Number of chargers,
(iv) Charger type,
(v) Dynamic charging demand over time in a day,
(vi) Partial recharging,
(vii) Queuing,
(viii) Multiple fixed routes,
MPCC: Mathematical programs with complementarity constraints,
N/A: Not applicable,
SA: Simulated annealing,
GA: Genetic algorithm,
CS: Commercial solver,
H: Heuristic,
NS: Not specified.}
\end{table}

\section{Problem definition}\label[sec]{methodology}
We now formally describe the EVSELCA problem, which  we model as an MILP. To ease reading,  we use calligraphic letters to represent sets (e.g., $\sR$), uppercase Roman letters for parameters (e.g., $\overline{T}$), lowercase Roman letters for variables and indices (e.g., $y_f$), and Greek letters (e.g., $\alpha$) as superscripts to modify parameters and variables. 

Let $\sC$ denote a set of customers, and let $\sR$ denote a set of trucks serving these customers with routes to be electrified while keeping their routes intact. (That is, regenerating routes from scratch due to electrification is not of interest.) We call these routes \textit{EV routes}. The goal of the EVSELCA is to allow EVs to complete their daily operations at a minimum cost by planning recharging infrastructure and scheduling recharging activities. Although solving a strategic decision-making problem, the model respects operational limitations, such as total route time and charging capacities. Each EV route $r\in\sR$ contains a subset of customers, denoted by the subset $\sC_r=\{c_{0r}, c_{1r}, ..., c_{ir}, ..., c_{Nr}\}$ satisfying $\forall i \neq j, (i, j) \neq (0, N)$ by the order of visits: $c_{0r}\rightarrow c_{1r}\rightarrow ..., c_{N-1,r}\rightarrow c_{Nr}$, where $c_{0r} = c_{Nr} =$ depot of route $r\in \sR$. A set of charger types $\sK$ can be located at a set of candidate facilities $\sF$. The optimization time is discretized into a set of time steps $\sT$. \cref{set} and \cref{param} provide sets and parameters used in the model, respectively. Note that the definition of sets will be modified in \cref{transformation} as we will cluster the customers to simplify the problem. 

We now state our critical modeling assumptions (and note that some of these can easily be relaxed). We assume customer demand is deterministic, which means that the travel sequences in the routes remain unchanged over time except for accommodating on-route charging activities. Furthermore, the adoption of EV technology does not change the truck routes, and the customer visit sequences of a route are assumed to be the same as in conventional truck routes. Additionally, EVs are assumed to be identical. The amount of energy received from a charger type is assumed to be a linear function of recharging time, and the energy spent by EVs is assumed to be a linear function of travel time. Furthermore, we assume there is a fixed set of locations for possible EV charging facilities.

\begin{table}[!htb]
  \footnotesize
  \caption{Sets used in the MILP.}\label[tab]{set}
  \begin{tabularx}{\linewidth}{lX}
  \toprule
    \textbf{Set} & \textbf{Definition}\\
    \midrule
    $\sC$ & set of customers, $\sC = \cup_{r\in \sR}{\sC_r}\replace{}{, \sC_r \cap \sC_{r^\prime} = \emptyset \quad \forall r, r^\prime \in \sR}$\\
    $\sC_r$ & subset of customers in route $r\in \sR$, $\sC_r = \{c_{0r}, c_{1r}, ..., c_{ir}, ..., c_{Nr}\}$ satisfy $\forall i \neq j, (i, j) \neq (0, N)$ by continuation order of visit: $c_{0r}\rightarrow c_{1r}\rightarrow ..., c_{N-1,r}\rightarrow c_{Nr}$, where $c_{0r} = c_{Nr} =$ depot of route $r\in \sR$\\
    $\sC_{rft}$ & subset of customers in route $r\in\sR$, after serving which the vehicle may visit $f\in\sF$ for recharging at time step $t\in\sT$, that is $\sC_{rft} = \Big\{c\in\sC_r | \sum_{j=0}^{i-1} \left[T_{c_{jr},c_{j+1,r}}^\tau + T_{c_{jr}}^\kappa\right] + T_{c_{ir}}^\kappa \leq D_t \leq \overline{T} - \sum_{j=i}^{N} \left[T_{c_{jr},c_{j+1,r}}^\tau + T_{c_{jr}}^\kappa\right]\Big\}  \forall r\in\sR,f\in\sF,t\in\sT$\\
    $\sF$ & set of candidate EV charging facility locations\\
    $\sF_{c_{ir}}$ & subset of candidate locations that the vehicle of $r\in\sR$, after serving customer $c_{ir}$, may possibly visit for recharging, that is, $\sF_{c_{ir}} = \Big\{f\in\sF | T_{c_{ir}f}^\tau \leq \overline{B}\land T_{c_{ir}f}^\tau + T_{fc_{i+1,r}}^\tau \leq \overline{T} - \sum_{j=0}^{N-1} \left[T_{c_{jr},c_{j+1,r}}^\tau + T_{c_{jr}}^\kappa\right] +\min\{0,{B_r^\iota-\sum_{j=0}^{N-1} {T_{c_{jr},c_{j+1,r}}^\tau}-B_r^\omega\over {R_K}}  \}  \Big\}  \forall c_{ir}\in \sC_r \setminus \{c_{Nr}\},r\in\sR \bigcup \Big\{f\in\sF | T_{c_{Nr}f}^\tau \leq \overline{B}\land T_{c_{Nr}f}^\tau + T_{fc_{Nr}}^\tau \leq \overline{T} - \sum_{j=0}^{N-1} \left[T_{c_{jr},c_{j+1,r}}^\tau + T_{c_{jr}}^\kappa\right] +\min\{0,{B_r^\iota-\sum_{j=0}^{N-1} {T_{c_{jr},c_{j+1,r}}^\tau}-B_r^\omega\over {R_K}}  \}  \Big\}  \forall r\in\sR$\\
    $\sK$ & set of charger types, $\sK = \{1,2,3,...,K\}$\\
    $\sR$ & set of routes\\
    $\sT$ & set of time steps\\
    $\sT_{c_{ir}f}$ & subset of time steps in the beginning of which the vehicle of $r\in\sR$ after serving customer $c_{ir}$ may possibly visit $f\in\sF_{c_{ir}}$ for recharging, that is, $\sT_{c_{ir}f} = \Big\{t\in\sT | \sum_{j=0}^{i-1} \left[T_{c_{jr},c_{j+1,r}}^\tau + T_{c_{jr}}^\kappa\right] + T_{c_{ir}}^\kappa \leq D_t \leq \overline{T} - \sum_{j=i}^{N-1} \left[T_{c_{jr},c_{j+1,r}}^\tau + T_{c_{jr}}^\kappa\right]\Big\}  \forall c_{ir}\in \sC_r,r\in\sR,f\in\sF_{c_{ir}}$\\
    \bottomrule
  \end{tabularx}
\end{table}

\begin{table}[!htbp]
  \footnotesize
  \caption{Parameters used in the MILP.}\label[tab]{param}
  \begin{tabularx}{\linewidth}{lX}
  \toprule
    \textbf{Parameter} & \textbf{Definition}\\
    \midrule
    $\overline{B}$ & maximum battery capacity in time unit\\
    $B_r^\iota$ & initial battery capacity of EV route $r\in\sR$\\
    $B_r^\omega$ & desired final battery capacity of EV route $r\in\sR$\\
    $C^\rho$ & value of time spent for recharging and driving to recharging facility\\
    $C^\xi_k$ & energy cost of recharging per time unit with charger type $k\in\sK$\\
    $C_k^\nu$ & cost of installing charger type $k\in\sK$ per time unit\\
    $C_f^\phi$ & fixed charging facility cost at candidate location $f\in\sF$ per time unit\\
    $D_t$ & actual time associated with time step $t\in\sT$, that is $D_t = tT^\Delta$ \\
    $\bigM$ & an adequately large number, e.g., $\bigM>2\overline{B}$\\
    $R_k$ & recharging amount received from charger type $k\in\sK$ per time unit\\
    $\overline{T}$ & maximum allowed operational time for each EV route\\
    $T^\Delta$ & duration of time steps\\
    $T_{c_{ir}c_{jr}}^\tau$ & travel time from customer $c_{ir}\in \sC_r$ to customer $c_{jr}\in \sC_r$ by EV route $r\in\sR$\\
    $T_{c_{ir}f}^\delta$ & detour travel time from customer $c_{ir}\in \sC_r$ by EV route $r\in\sR$ to $f\in \sF$\\
    $T_{c_{ir}}^\kappa$ & time spent for serving customer $c_{ir}\in \sC_r$ by EV route $r\in \sR$\\
    \bottomrule
  \end{tabularx}
\end{table}

We let the binary variable $x_{c_{ir}fkt}=1$ denote that the EV of route $r \in \sR$ recharges at $f \in \sF$ using charger type $k \in \sK$ at time step $t \in \sT$ immediately after serving customer $c_{ir} \in \sC_r$; $x_{c_{ir}fkt}=0$, otherwise. To ensure $x_{c_{ir}fkt}$ represents the desired value, we keep track of the start and end time for recharging an EV for $c_{ir}$, $f$, and $k$ with continuous variables $s_{c_{ir}fk}$ and $e_{c_{ir}fk}$. Let the binary variable $x_{c_{ir}fkt}^\alpha=1$ indicate that the EV route $r\in\sR$ is recharging at facility $f\in \sF$ using charger $k\in\sK$ after detouring from $c_{ir} \in \sC_r$ for time steps, where the associated time $D_t$ falls within the range $[0, s_{c_{ir}fk})$ (excluding $t \in \sT$). On the other hand, let binary variable $x_{c_{ir}fkt}^\beta=1$ denote that the EV route $r\in\sR$ is recharging at facility $f\in \sF$ using charger $k\in\sK$ after detouring from $c_{ir} \in \sC_r$ for time steps, where the associated time $D_t$ falls within the range $[0, e_{c_{ir}fk}]$ (including $t \in \sT$). 

\begin{figure}[t]
   
  \hspace{0.5 cm} \includegraphics[width=.7\linewidth]{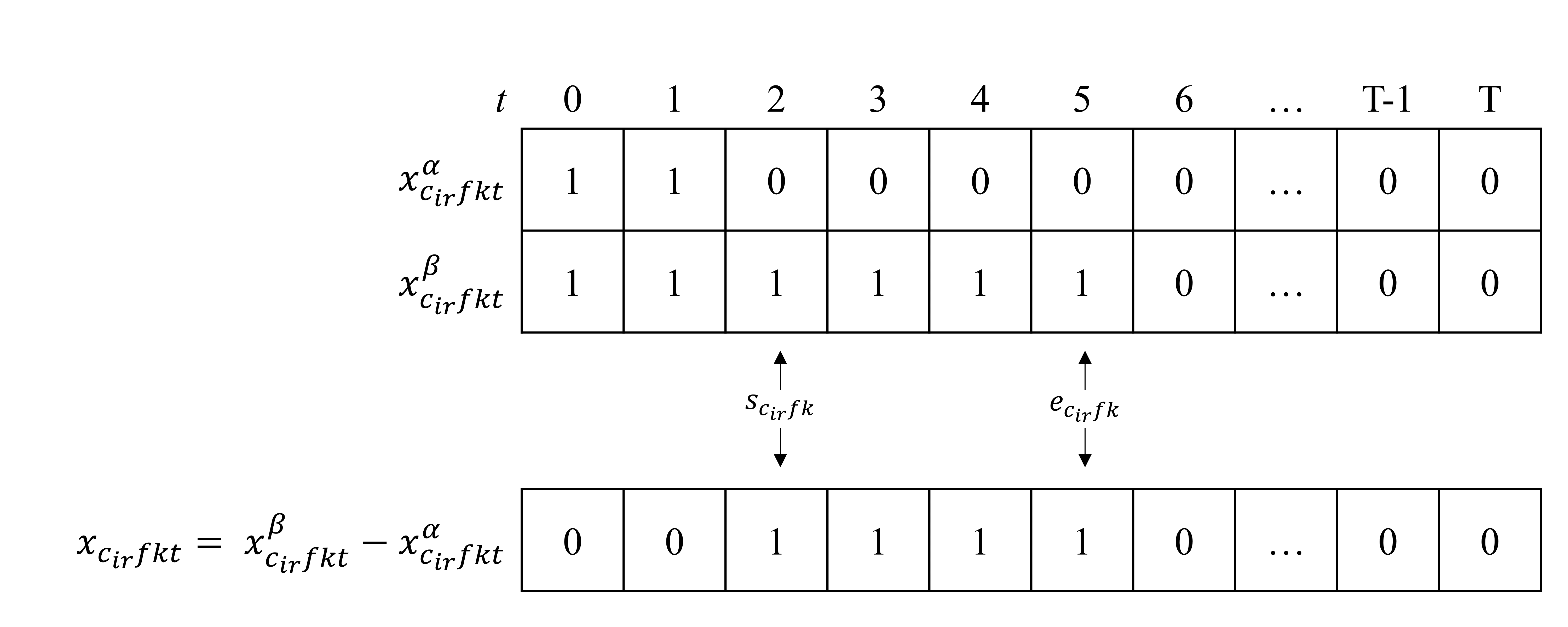}
  \caption{Example of how $x_{c_{ir}fkt}$ is maintained.\label[fig]{x_formulation_fig}}
\end{figure}

\cref{x_formulation_fig} illustrates how our constraints ensure the variable $x_{c_{ir}fkt}$ takes the desired value. In this example, EV route $r\in\sR$ starts recharging at facility $f\in \sF$ using charger $k\in\sK$ after detouring from $c_{ir} \in \sC_r$ at time $s_{c_{ir}fk}$, where the equivalent time step is $2$ and ends at time $e_{c_{ir}fk}$, where the equivalent time step is $5$\replace{}{, respectively}. The binary variable $x_{c_{ir}fkt}^\alpha=1$ for all time steps preceding and \textit{excluding} time step $2$, while the binary variable $x_{c_{ir}fkt}^\beta=1$ for all time steps preceding and \textit{including} time step $5$\replace{}{, respectively}. Consequently, utilizing the equation $x_{c_{ir}fkt} = x_{c_{ir}fkt}^{\beta} - x_{c_{ir}fkt}^{\alpha}$, the binary variable $x_{c_{ir}fkt}$ is 1 only for the expected time steps that are in the interval $[2, 5]$. 

The binary variable $y_f=1$ denotes $f \in \sF$ is open; $y_f=0$, otherwise. The variable $z_{fk}\in\mathbb{Z}_{\geq 0}$ denotes the number of charger type $k \in \sK$ allocated to $f \in \sF$.  Refer to \cref{vars} for the definitions of all variables used in the model. The EVSELCA problem is formulated as follows:    

\begin{table}[!htb]
  \footnotesize
  \caption{Variables used in the MILP.}\label[tab]{vars}
  \begin{tabularx}{\linewidth}{lX}
  \toprule
    \textbf{Variable} & \textbf{Definition}\\
    \midrule
    $b_{c_{ir}f}^\prime$ & remaining battery capacity (in time units) before arriving $f\in\sF$ for EV route $r\in\sR$ after detouring from $c_{ir} \in \sC_r$\\
    $b_{c_{ir}}$ & remaining battery capacity (in time units) before serving $c_{ir} \in \sC_r$ of EV route $r\in\sR$\\
    $d_{c_{ir}}$ & departure time from $c_{ir} \in \sC_r$ of EV route $r\in\sR$\\
    $u_{c_{ir}fk}$ & duration of recharging time EV route $r\in\sR$ spends at $f\in \sF$ using $k\in \sK$ after detouring from $c_{ir} \in \sC_r$\\
    $w_{c_{ir}fk}$ & duration of waiting time EV route $r\in\sR$ spends to recharge at $f\in \sF$ using $k\in \sK$ after detouring from $c_{ir} \in \sC_r$\\
    $s_{c_{ir}fk}$ & time that EV route $r\in\sR$ starts recharging at $f\in \sF$ using $k\in \sK$ after detouring from $c_{ir} \in \sC_r$\\
    $e_{c_{ir}fk}$ & time that EV route $r\in\sR$ ends recharging $f\in \sF$ using $k\in \sK$ after detouring from $c_{ir} \in \sC_r$\\
    $q_{c_{ir}fk}$ & $\begin{cases}
          1 & \text{\parbox[t][][t]{0.80\textwidth}{if EV route $r\in\sR$ detours from $c_{ir} \in \sC_r$ to $f\in \sF$ to recharge using $k\in\sK$}}\\
          0 & \text{otherwise}\\
         \end{cases}$\\
    $x_{c_{ir}fkt}^{\alpha}$ & $\begin{cases}
          1 & \text{\parbox[t][][t]{0.80\textwidth}{if EV route $r\in\sR$ is recharging at facility $f\in \sF$ using charger $k\in\sK$ after detouring from $c_{ir} \in \sC_r$ for time steps where the equivalent time, $D_t$, falls within the range $[0, s_{c_{ir}fk})$ (excluding $t \in \sT$)}} \\
          0 & \text{otherwise}\\
         \end{cases}$\\
         
    $x_{c_{ir}fkt}^{\beta}$ & $\begin{cases}
          1 & \text{\parbox[t][][t]{0.80\textwidth}{if EV route $r\in\sR$ is recharging at facility $f\in \sF$ using charger $k\in\sK$ after detouring from $c_{ir} \in \sC_r$ for time steps where the equivalent time, $D_t$, falls within the range $[0, e_{c_{ir}fk}]$ (including $t \in \sT$)}} \\
          0 & \text{otherwise}\\
         \end{cases}$\\
    $x_{c_{ir}fkt}$ & $\begin{cases}
          1 & \text{\parbox[t][][t]{0.80\textwidth}{if EV route $r\in\sR$ is recharging at  $f\in \sF$ using $k\in\sK$ at time step $t\in \sT$ after detouring from $c_{ir} \in \sC_r$}} \\
          0 & \text{otherwise}\\
         \end{cases}$\\
    $y_f$ & $\begin{cases}
          1 & \text{if a charging facility is located at $f\in\sF$} \\
          0 & \text{otherwise}\\
         \end{cases}$\\
    $z_{fk}$ & number of charger type $k\in\sK$ installed at $f\in\sF$, $z_{fk}\in\mathbb{Z}_{\geq 0}$\\
    \bottomrule
  \end{tabularx}
\end{table}
\begin{equation}\label[objfunc]{obj_fun}
\begin{split}
    & \begin{multlined}
        \min  \textbf{C} = \sum_{\substack{c_{ir}\in \sC_r,r \in \sR,\\f\in \sF_{c_{ir}}, k\in \sK}} \left[C^\rho\left(T_{c_{ir}f}^\delta q_{c_{ir}fk}+w_{c_{ir}fk}\right)+(C^\rho +C^\xi_k) u_{c_{ir}fk}\right] + \sum_{f\in\sF} C_f^\phi y_f \\ + \sum_{f\in\sF,k\in\sK}C_k^\nu z_{fk}
    \end{multlined}
\end{split}
\end{equation}
\noindent subject to,
\begin{equation}\label[consset]{x_formulation}
    x_{c_{ir}fkt} = x_{c_{ir}fkt}^{\beta} - x_{c_{ir}fkt}^{\alpha} \qquad \forall {c_{ir}}\in\sC_r, r\in \sR, f\in \sF_{c_{ir}}, k\in \sK, t\in \sT_{c_{ir}f}
\end{equation}
\begin{equation}\label[consset]{charger_availability}
    \sum_{{c_{ir}}\in\sC_{rft}, r\in\sR} x_{c_{ir}fkt} \leq z_{fk} \qquad \forall f\in \sF, k\in \sK, t\in \sT
\end{equation}
\begin{equation}\label[cons]{u_time_step}
    u_{c_{ir}fk}\leq T^\Delta \sum_{t\in\sT_{c_{ir}f}} x_{c_{ir}fkt} \qquad \forall {c_{ir}}\in\sC_r, r\in \sR, f\in \sF_{c_{ir}}, k\in \sK
\end{equation}
\begin{equation}\label[cons]{u_time_step2}
    T^\Delta \sum_{t\in\sT_{c_{ir}f}} x_{c_{ir}fkt} - u_{c_{ir}fk} \leq T^\Delta - \epsilon \qquad \forall {c_{ir}}\in\sC_r, r\in \sR, f\in \sF_{c_{ir}}, k\in \sK
\end{equation}
\begin{equation}\label[cons]{have_chargers_if_open}
    z_{fk} \leq \bigM y_f \qquad \forall f\in \sF, k\in \sK
\end{equation}
\begin{equation}\label[cons]{charge_u_if_decided_to_charge}
    u_{c_{ir}fk}\leq \bigM q_{c_{ir}fk} \qquad \forall {c_{ir}}\in\sC_r, r\in \sR, f\in \sF_{c_{ir}}, k\in \sK
\end{equation}
\begin{equation}\label[cons]{wait_w_if_decided_to_charge}
    w_{c_{ir}fk} \leq \bigM q_{c_{ir}fk} \qquad \forall {c_{ir}}\in\sC_r, r\in \sR, f\in \sF_{c_{ir}}, k\in \sK
\end{equation}
\begin{equation}\label[cons]{detour_if_decided_to_charge1}
    q_{c_{ir}fk} \leq \sum_{t\in\sT_{c_{ir}f}}x_{c_{ir}fkt} \qquad \forall {c_{ir}}\in\sC_r, r\in \sR, f\in \sF_{c_{ir}}, k\in \sK
\end{equation}
\begin{equation}\label[cons]{detour_if_decided_to_charge2}
    x_{c_{ir}fkt} \leq q_{c_{ir}fk} \qquad \forall {c_{ir}}\in\sC_r, r\in \sR, f\in \sF_{c_{ir}}, k\in \sK, t\in\sT_{c_{ir}f}
\end{equation}
\begin{equation}\label[cons]{charge_in_one_charger_and_one_station}
    \sum_{ f\in\sF_{c_{ir}}, k\in\sK}{q_{c_{ir}fk}} \leq 1 \qquad \forall {c_{ir}}\in\sC_r, r\in \sR
\end{equation}
\begin{equation}\label[consset]{battery_capacity_flow1}
    b^\prime_{c_{ir}f} \leq b_{c_{ir}}-T_{c_{ir}f}^\tau +\bigM (1-\sum_{k\in\sK}{q_{c_{ir}fk}}) \qquad \forall {c_{ir}}\in\sC_r, r\in \sR, f\in\sF_{c_{ir}}
\end{equation}
\begin{equation}\label[consset]{battery_capacity_flow2}
    b^\prime_{c_{ir}f} \leq \bigM \sum_{k\in\sK}{q_{c_{ir}fk}} \qquad \forall {c_{ir}}\in\sC_r, r\in \sR, f\in\sF_{c_{ir}}
\end{equation}
\begin{equation}\label[consset]{battery_capacity_flow2_2}
    b^\prime_{c_{ir}f} + \sum_{k\in\sK} u_{c_{ir}fk} \leq \overline{B} \qquad \forall {c_{ir}}\in\sC_r, r\in \sR, f\in\sF_{c_{ir}}
\end{equation}
\begin{equation}\label[consset]{battery_capacity_flow3}
\begin{split}
    & \begin{multlined}
        b_{c_{i+1,r}} \leq b_{c_{ir}}-T_{c_{ir},c_{i+1,r}}^\tau(1-\sum_{f\in\sF_{c_{ir}},k\in \sK}{q_{c_{ir}fk}})+\bigM \sum_{f\in\sF_{c_{ir}},k\in \sK}q_{c_{ir}fk} \\ \qquad \forall {c_{ir}}\in\sC_r, r\in \sR
    \end{multlined}
\end{split}
\end{equation}
\begin{equation}\label[consset]{battery_capacity_flow4}
\begin{split}
    & \begin{multlined}
    b_{c_{i+1,r}} \leq \sum_{f\in\sF_{c_{ir}}} \left[b^\prime_{c_{ir}f}+\sum_{k\in\sK}\left(R_ku_{c_{ir}fk} - {T_{fc_{i+1,r}}^\tau q_{c_{ir}fk}}\right)\right] \\ +\bigM \Big(1-\sum_{f\in\sF_{c_{ir}},k\in\sK}{q_{c_{ir}fk}}\Big) \quad \forall {c_{ir}}\in\sC_r, r\in \sR
    \end{multlined}
\end{split}
\end{equation}
\begin{equation}\label[consset]{battery_initial_capacity}
    b_{c_{0r}} = B_r^\iota \qquad \forall {c_{ir}}\in\sC_r, r\in \sR
\end{equation}
\begin{equation}\label[consset]{battery_final_capacity}
    b_{c_{Nr}} \geq B_r^\omega \qquad \forall {c_{ir}}\in\sC_r, r\in \sR
\end{equation}
\begin{equation}\label[consset]{time}
\begin{split}
& 
\begin{multlined} d_{c_{ir}} = \sum_{j=0}^{i-1} \left[T_{c_{jr},c_{j+1,r}}^\tau+\sum_{f\in\sF_{c_{ir}},k\in\sK}(T_{c_{jr}f}^\delta q_{c_{jr}fk} + u_{c_{jr}fk} + w_{c_{jr}fk})+ T_{c_{jr}}^\kappa\right] + T_{c_{ir}}^\kappa  \\ \forall {c_{ir}}\in\sC_r, r\in \sR
\end{multlined}
\end{split}
\end{equation}
\begin{equation}\label[consset]{charging_start_time_formulation}
\begin{split} & 
    \begin{multlined}
    s_{c_{ir}fk} = d_{c_{ir}} + T_{c_{ir}f}^\tau {q_{c_{ir}fk}} + {w_{c_{ir}fk}}\\ \forall {c_{ir}}\in\sC_r, r\in \sR, f\in\sF_{c_{ir}}, k\in \sK, t\in \sT_{c_{ir}f}
    \end{multlined}
\end{split}
\end{equation}
\begin{equation}\label[consset]{charging_end_time_formulation}
\begin{split} & 
    \begin{multlined}
    e_{c_{ir}fk} = d_{c_{ir}} + T_{c_{ir}f}^\tau{q_{c_{ir}fk}} + {w_{c_{ir}fk}} + {u_{c_{ir}fk}}\\ \forall {c_{ir}}\in\sC_r, r\in \sR, f\in\sF_{c_{ir}}, k\in \sK, t\in \sT_{c_{ir}f}
    \end{multlined}
\end{split}
\end{equation}
\begin{equation}\label[consset]{x_alpha_formulation1}
\begin{split}
    & \begin{multlined}
     s_{c_{ir}fk} \leq D_t + T^\Delta-\epsilon + \bigM \left(1-{q_{c_{ir}fk}}+x_{{c_{ir}fkt}}^{\alpha}\right)\\ \forall {c_{ir}}\in\sC_r, r\in \sR, f\in\sF_{c_{ir}}, k\in \sK, t\in \sT_{c_{ir}f}
    \end{multlined}
\end{split}
\end{equation}
\begin{equation}\label[consset]{x_alpha_formulation2}
\begin{split}
    & \begin{multlined}
     s_{c_{ir}fk} \geq D_t + T^\Delta- \bigM (2-{q_{c_{ir}fk}}-x_{{c_{ir}fkt}}^{\alpha})\\ \forall {c_{ir}}\in\sC_r, r\in \sR, f\in\sF_{c_{ir}}, k\in \sK, t\in \sT_{c_{ir}f}
    \end{multlined}
\end{split}
\end{equation}
\begin{equation}\label[consset]{x_end_formulation1}
\begin{split}
    & \begin{multlined}
     e_{c_{ir}fk} \leq D_t -\epsilon + \bigM \left(1-{q_{c_{ir}fk}}+x_{{c_{ir}fkt}}^{\beta}\right)\\ \qquad \forall {c_{ir}}\in\sC_r, r\in \sR, f\in\sF_{c_{ir}}, k\in \sK, t\in \sT_{c_{ir}f}
    \end{multlined}
\end{split}
\end{equation}
\begin{equation}\label[consset]{x_end_formulation2}
\begin{split}
    & \begin{multlined}
         e_{c_{ir}fk} \geq D_t - \bigM \left(2-{q_{c_{ir}fk}}-x_{{c_{ir}fkt}}^{\beta}\right)\\ \qquad \forall {c_{ir}}\in\sC_r, r\in \sR, f\in\sF_{c_{ir}}, k\in \sK, t\in \sT_{c_{ir}f}
    \end{multlined}
\end{split}
\end{equation}
\begin{equation}\label[cons]{operation_time_constraint}
    d_{c_{Nr}} \leq \overline{T} \qquad \forall r\in \sR
\end{equation}
\begin{equation*}\label[consset]{non-neg}
    \nonumber x_{c_{ir}fkt}, ~x_{c_{ir}fkt}^{\alpha}, ~x_{c_{ir}fkt}^{\beta} ,~y_f, ~q_{c_{ir}fk} \in \{0,1\}, z_{fk}\in\mathbb{Z}_{\geq 0}, b_{c_{ir}}, ~d_{c_{ir}}, ~u_{c_{ir}fk}, ~w_{c_{ir}fk}\in \mathbb{R}_{\geq 0} .
\end{equation*}
\smallskip

The objective function \labelcref{obj_fun} minimizes the costs associated with detouring, waiting, recharging, energy consumption, charging facility, and charger installation. Note that all costs are minimized per time unit (e.g., one day). Constraint \labelcref{x_formulation} defines the variable $x_{c_{ir}fkt}$ in terms of $x_{c_{ir}fkt}^\alpha$ and $x_{c_{ir}fkt}^\beta$. Constraint \labelcref{charger_availability} ensures that the total number of trucks charging at a time step does not exceed the capacity of the charging facility. Constraint \labelcref{u_time_step} satisfies that the recharging time does not exceed the total time occupied by a truck at a facility, and constraint \labelcref{u_time_step2} ensures that a truck does not occupy a charger while not recharging. Constraint \labelcref{have_chargers_if_open} guarantees that if a charging facility is not located at $f$, a charger should not be allocated. Constraints \labelcref{charge_u_if_decided_to_charge} and \labelcref{wait_w_if_decided_to_charge} enforce the charging time $u_{c_{ir}fk}$ and waiting time $w_{c_{ir}fk}$ to be zero when recharging does not occur. Constraints \labelcref{detour_if_decided_to_charge1} and \labelcref{detour_if_decided_to_charge2} ensure that if the truck recharges, the charger will be considered occupied at least in one time step. Constraint \labelcref{charge_in_one_charger_and_one_station} ensures that a truck can recharge only at one facility using one charger type after serving a customer. 

Constraints \labelcref{battery_capacity_flow1}–\labelcref{battery_capacity_flow4} define the remaining battery of trucks (in time units) after serving every customer ($b_{c_{ir}}$) and immediately before arriving at a charging facility ($b^\prime_{c_{ir}f}$). Constraints \labelcref{battery_initial_capacity} and \labelcref{battery_final_capacity} enforce the initial and desired final battery capacity to be equal to $B_r^{\iota}$ and $B_r^{\omega}$, respectively. Constraint \labelcref{time} defines the departure time of trucks after serving each customer. Constraints \labelcref{x_alpha_formulation1} and \labelcref{x_alpha_formulation2} define the variable $x_{{c_{ir}fkt}}^{\alpha}$. Similarly, constraints \labelcref{x_end_formulation1} and \labelcref{x_end_formulation2} define the variable $x_{{c_{ir}fkt}}^{\beta}$. Constraint \labelcref{operation_time_constraint} ensures that the total operational time cannot exceed the maximum allowed operational time for each EV route. 

\section{Metaheuristic solution approaches with clustering}\label[sec]{heuristic}

We now present a metaheuristic solution procedure for solving the EVSELCA problem, as defined by Equations \labelcref{obj_fun}–-\labelcref{operation_time_constraint}. 
We employ a metaheuristic approach because of substantial growth in solution space of the EVSELCA problem as the problem size increases.
First, we adopt a clustering approach in \cref{clustering} to simplify the problem. This involves redefining sets, parameters, and variables in the MILP so that it can be used with clusters, as shown in \cref{transformation}. Although clustering can help address the computational difficulty at a smaller scale, we further develop a metaheuristic solution method using the GA to tackle the EVSELCA problem in \cref{GA_section}. Then, in \cref{hybrid_sec} we describe a hybrid solution approach that combines the GA and MILP solvers.

\subsection{Clustering}\label[sec]{clustering}
EVSELCA problems (Equations \labelcref{obj_fun}–\labelcref{operation_time_constraint}) can be complex, making solving them impractical for large-scale instances. To overcome this challenge, we adopt a clustering approach proposed by \cite{cokyasar2023}. This approach simplifies the problem by grouping customers into clusters and limiting recharging to after the completion of service at these clusters. That is, rather than considering recharging after any customer, we create clusters of customers and assume recharging occurs only after completing service at these clusters.

The clustering method aims to identify the best cut-points with a given number of clusters to satisfy the following conditions:

\begin{enumerate}[label=\roman*.]
  \item All customers in a cluster must belong to a single route. 

  \item The intersection of any two clusters must be empty.
  \item Customers in a cluster must follow the order of service.
  \item The distance traveled within a cluster must not exceed a certain threshold, such as a portion of the EV range.
\end{enumerate}

 This clustering method uses an optimization model that maximizes the spatial difference between clusters. In this model, the binary variable, $p_{nc_{ir}} = 1$ indicates cut $n\in\sN$ is placed right after customer $c_{ir}$. \replace{The auxiliary binary variable $m_{nc_{ir}}$ regulates the order of cuts and conserves sequencing.}{The auxiliary binary variable $m_{nc_{ir}} = 1$ for all customers belonging to $\{c_{1r}, c_{2r}, \dots, c_{ir}\}$ if cut $n\in\sN$ is placed right after customer $c_{ir}\in\sC_r$, and $m_{nc_{ir}} = 0$ for all customers belonging to $\{c_{i+1,r}, c_{i+2,r}, \dots, c_{N-1,r}\}$.} \cref{clustering_set_var} provides sets and parameters used in the clustering model. The clustering optimization model is formulated as follows. 

\begin{equation}\label[objfunc]{clustering_obj_fun}
    \max \sum_{\substack{c_{ir}\in \sC_r,r \in \sR,\\n\in \sN}} T_{c_{ir}c_{jr}}^\tau p_{nc_{ir}} 
\end{equation}
\noindent subject to,
\begin{equation}\label[consset]{p_formulation}
    p_{nc_{ir}} = m_{nc_{ir}} - m_{nc_{i+1,r}} \qquad \forall {c_{ir}}\in\sC_r, r\in \sR, n\in \sN
\end{equation}
\begin{equation}\label[consset]{m_formulation}
    \replace{m_{nc_{i+1,r}}}{m_{n+1,c_{ir}}} - m_{nc_{ir}} \geq 0 \qquad \forall {c_{ir}}\in\sC_r, r\in \sR, n\in \sN
\end{equation}
\begin{equation}\label[consset]{p_sumoverc}
    \replace{\sum_{\substack{c_{ir}\in \sC_r,r \in \sR}}p_{nc_{ir}} = 1 \qquad \forall n\in \sN}{\sum_{\substack{c_{ir}\in \sC_r}}p_{nc_{ir}} = 1 \qquad \forall r \in \sR, n\in \sN}
\end{equation}
\begin{equation}\label[consset]{p_sumovern}
    \sum_{n\in \sN}p_{nc_{ir}} \leq 1 \qquad \forall {c_{ir}}\in\sC_r, r\in \sR
\end{equation}

The objective function \labelcref{clustering_obj_fun} maximizes the total travel time between clusters to ensure that clusters are sufficiently apart from each other. Constraints \labelcref{p_formulation} and \labelcref{m_formulation} conserve the order of cuts, preventing cut $n$ from being placed after cut $n+1$. Constraint \labelcref{p_sumoverc} ensures that every cut is positioned immediately after one specific point, while constraint \labelcref{p_sumovern} ensures that only one cut can be placed after a particular point. The model in \citet{cokyasar2023} also presents the within the cluster travel time constraints. As these constraints require more parametric definition, we refrain providing these constraints for simplicity.

With a predetermined number of clusters, the clustering method may not always guarantee a feasible solution because the travel time within clusters could exceed a preset threshold. In our analysis, we start with a small $|\sN|$, solve the problem, and increment $|\sN|$ by one until a feasible solution is obtained. Therefore, we find the minimum number of clusters ($|\sN|> 0$) and their partitioning. Note that such an approach sacrifices the solution quality to gain solution speed. Our analysis in the upcoming sections will be a product of this sacrifice, and a real case should better be handled with $|\sN|$ large enough to obtain a higher quality solution.

\begin{table}[!htb]
  \footnotesize
  \caption{Sets and variables used in the clustering model.}\label[tab]{clustering_set_var}
  \begin{tabularx}{\linewidth}{lX}
  \toprule
    \textbf{Set} & \textbf{Definition}\\
    \midrule
    $\sN$ & set of cut-points for clustering \\
    \midrule
    \textbf{Variable} & \textbf{Definition}\\
    \midrule
        $p_{nc_{ir}}$ & $\begin{cases}
        1 & \text{\parbox[t][][t]{0.80\textwidth}{if  cut $n \in \sN$ is placed right after customer $c_{ir}$}}\\
        0 & \text{otherwise}\\
        \end{cases}$\\
        $m_{nc_{ir}}$ & auxiliary binary variable regulating the order of cuts \\
    \bottomrule
\end{tabularx}
\end{table}

\subsection{Transformation}\label[sec]{transformation}
In switching from customer to cluster, the definition of some sets, parameters, variables, and constraints changes. It is straightforward to redefine $\sC$ from being the set of individual customers to being the set of clusters of customers. We then define a new parameter $T_{c_{ir}}^\gamma$ to represent the total travel time  based on the order of visits within the cluster for serving customers of cluster $c_{ir}$ by the EV route $r$. Moreover, we redefine parameters $T_{c_{ir}c_{jr}}^\tau$ to be the travel time from the last customer of cluster $c_{ir}$ to the first customer of cluster $c_{jr}$ and $T_{c_{ir}f}^\delta$ to be the travel time from the last customer of cluster $c_{ir}$ to the location of facility $f$. The model in \cref{methodology} could be structured by considering the clustering approach from the beginning, yet we prefer the given presentation to provide a general model without dependence on the clustering heuristic. 

With these changes, we redefine the subset $\sC_{rft}$ in \cref{Crft}. 
\begin{definition}\label[def]{Crft}
For a given $r\in\sR$, $f\in\sF$, and $t\in\sT$, 

\begin{equation}
    \nonumber
    \begin{split}
        & \begin{multlined}
        \sC_{rft} = \Big\{c\in\sC_r | \sum_{j=0}^{i-1} \left[T_{c_{jr},c_{j+1,r}}^\tau + T_{c_{jr}}^\kappa + T_{c_{jr}}^\gamma\right] + T_{c_{ir}}^\kappa + T_{c_{ir}}^\gamma \\ \leq t \leq \overline{T} - \sum_{j=i}^{N} \left[T_{c_{jr},c_{j+1,r}}^\tau + T_{c_{jr}}^\kappa + T_{c_{ir}}^\gamma\right]\Big\}.  
        \end{multlined}
    \end{split}
\end{equation}
\end{definition}
The definition of $\sC_{rft}$ requires that an EV route can  recharge only after serving the last customer of the cluster, rather than in the middle of the cluster. Therefore, $\sC_{rft}$ refers to the subset of clusters on route $r\in\sR$ that the vehicle must complete before it can visit $f\in\sF$ for recharging at time step $t\in\sT$. 

Furthermore, we redefine $\sF_{c_{ir}}$ and $\sT_{c_{ir}f}$ due to the changes made to the definition of parameters. The updated $\sF_{c_{ir}}$ denotes the subset of facilities that the vehicle of route $r\in\sR$ can visit after serving customers of cluster $c_{ir}$, as defined in \cref{Fcir}. Similarly, $\sT_{c_{ir}f}$ denotes the subset of time steps during which the vehicle of route $r\in\sR$ can visit $f\in\sF_{c_{ir}}$ for recharging after serving customers of cluster $c_{ir}$, as defined in \cref{Tcirf}.   

\begin{definition}\label[def]{Fcir}
For a given $c_{ir}\in \sC_r \setminus \{c_{Nr}\}$ and $r\in\sR $, 

\begin{equation}
\nonumber
    \begin{split}
        & \begin{multlined}
        \sF_{c_{ir}} = \Big\{f\in\sF | T_{c_{ir}f}^\tau \leq \overline{B}\land T_{c_{ir}f}^\tau + T_{fc_{i+1,r}}^\tau \leq \overline{T} - \sum_{j=0}^{N-1} \left[T_{c_{jr},c_{j+1,r}}^\tau + T_{c_{jr}}^\kappa + T_{c_{jr}}^\gamma\right] \\ +\min\{0,{B_r^\iota-\sum_{j=0}^{N-1} [{T_{c_{jr},c_{j+1,r}}^\tau + T_{c_{jr}}^\gamma}]-B_r^\omega\over {R_K}}  \}  \Big\}  \\ \bigcup \Big\{f\in\sF | T_{c_{Nr}f}^\tau \leq \overline{B}\land T_{c_{Nr}f}^\tau + T_{fc_{Nr}}^\tau \leq \overline{T} - \sum_{j=0}^{N-1} \left[T_{c_{jr},c_{j+1,r}}^\tau + T_{c_{jr}}^\kappa + T_{c_{jr}}^\gamma\right] \\ +\min\{0,{B_r^\iota-\sum_{j=0}^{N-1} [{T_{c_{jr},c_{j+1,r}}^\tau + T_{c_{jr}}^\gamma}]-B_r^\omega\over {R_K}}  \}  \Big\}.
        \end{multlined}
    \end{split}
\end{equation}
\end{definition}
\begin{definition}\label[def]{Tcirf}
For a given $c_{ir}\in \sC_r$, $r\in\sR$, and $f\in\sF_{c_{ir}}$, 

\begin{equation}
\nonumber
    \begin{split}
        & \begin{multlined}
           \sT_{c_{ir}f} = \Big\{t\in\sT | \sum_{j=0}^{i-1} \left[T_{c_{jr},c_{j+1,r}}^\tau + T_{c_{jr}}^\kappa 
+ T_{c_{jr}}^\gamma\right] + T_{c_{ir}}^\kappa + T_{c_{ir}}^\gamma \\ \leq t \leq \overline{T} - \sum_{j=i}^{N-1} \left[T_{c_{jr},c_{j+1,r}}^\tau + T_{c_{jr}}^\kappa + T_{c_{jr}}^\gamma\right]\Big\}.
        \end{multlined}
    \end{split}
\end{equation}
\end{definition}

Next, we replace the term \textit{customer} with \textit{cluster} wherever it is used in \cref{set}, \cref{param}, and \cref{vars}. Following the changes we have made to the sets, parameters, and variable definitions, we replace the constraints \labelcref{battery_capacity_flow1}, \labelcref{battery_capacity_flow3}, and \labelcref{time} with \labelcref{battery_capacity_flow1_modified}, \labelcref{battery_capacity_flow3_modified}, and \labelcref{time_modified}, respectively.

\begin{equation}\label[consset]{battery_capacity_flow1_modified}
    b^\prime_{c_{ir}f} \leq b_{c_{ir}}-T_{c_{ir}}^\gamma-T_{c_{ir}f}^\tau +\bigM (1-\sum_{k\in\sK}{q_{c_{ir}fk}}) \qquad \forall {c_{ir}}\in\sC_r, r\in \sR, f\in\sF_{c_{ir}}
\end{equation}
\begin{equation}\label[consset]{battery_capacity_flow3_modified}
\begin{split}
    & \begin{multlined}
    b_{c_{i+1,r}} \leq b_{c_{ir}}-(T_{c_{ir},c_{i+1,r}}^\tau+T_{c_{ir}}^\gamma)(1-\sum_{f\in\sF_{c_{ir}},k\in \sK}{q_{c_{ir}fk}})+\bigM \sum_{f\in\sF_{c_{ir}},k\in \sK}q_{c_{ir}fk} \\ \qquad \forall {c_{ir}}\in\sC_r, r\in \sR
    \end{multlined}
\end{split}
\end{equation}
\begin{equation}\label[consset]{time_modified}
\begin{split}
& 
\begin{multlined} d_{c_{ir}} = \sum_{j=0}^{i-1} \left[T_{c_{jr},c_{j+1,r}}^\tau+T_{c_{jr}}^\gamma+\sum_{f\in\sF_{c_{ir}},k\in\sK}(T_{c_{jr}f}^\delta q_{c_{jr}fk} + u_{c_{jr}fk} + w_{c_{jr}fk})+ T_{c_{jr}}^\kappa\right] \\ + T_{c_{ir}}^\kappa +T_{c_{ir}}^\gamma \qquad \forall {c_{ir}}\in\sC_r, r\in \sR
\end{multlined}
\end{split}
\end{equation}

\subsection{The genetic algorithm}\label[sec]{GA_section}
In this study we employ a tailored genetic algorithm  to solve the EVSELCA problem.
The GA is an evolutionary optimization search technique that has been widely used to solve MILPs \citep{Katoch2021}. \cref{GA_main} provides pseudocode for the approach; the functions used therein are detailed in the \hyperref[gaappendix]{Supplementary Material}. Key decision variables in the model are $q_{c_{ir}fk}$, $x_{c_{ir}fkt}$, $y_f$, and $z_{fk}$, which are the same as in \cref{vars} with the redefined $\sC$. Among key variables, $y_f$ and $z_{fk}$ relate to strategic decision-making, while $x_{c_{ir}fkt}$ aids in making an operational decision. At the tactical level, $q_{c_{ir}fk}$ plays an important role since it determines after which cluster to recharge, where to recharge, and what type of charger to use. To this end, we begin with exploring a solution for the tactical variable  that implicitly impacts solutions to strategic variables and provides implied time bounds for the recharging time. That is, a solution to other variables can be derived for given solutions to $q_{c_{ir}fk}$. First, $N^{pop}$ number of solutions for $q_{c_{ir}fk}$ is generated via the \textproc{Initialization} function as an initial population. In the initialization step we randomly select a number of clusters, following which a recharging is planned; and we select a facility for the recharging using a roulette wheel selection method (i.e., closer facilities have a higher chance of being selected). Once the \textit{where} aspect of $q_{c_{ir}fk}$ is addressed, we randomly select a type of charger for those facilities that were just picked to be visited. 
This population is then passed into \textproc{Crossover} and \textproc{Mutation} functions to potentially find a better solution. 

In \cref{GA_main}, $T_r^\rho$ was used to represent the route travel time minus the time spent serving customers, recharging, and waiting, as defined in \cref{Trrho}.

\begin{definition}\label[def]{Trrho}
For a given $r\in \sR$, 
\begin{equation}
    \nonumber
    T_r^\rho = \sum_{j=0}^{N-1} \left[T_{c_{jr},c_{j+1,r}}^\tau+\sum_{f\in\sF_{c_{ir}},k\in\sK}T_{c_{jr}f}^\delta q_{c_{jr}fk}\right].
\end{equation} 
\end{definition}

The value of $z_{fk}$ is estimated to calculate the objective value for each given $q_{c_{ir}fk}$ in the initial population. For a given $q_{c_{ir}fk}$ equal to $1$, $z_{fk}$ can get values between $1$ and $\sum_{c_{ir}\in C_{rft}, r\in \sR}{q_{c_{ir}fk}}$ as \cref{lemma2} denotes. (See the~\hyperref[proofs]{Appendix} for proofs of lemmas). The trade-off between the waiting time and facility cost depends on $z_{fk}$. The maximum of $z_{fk}$ ($\sum_{c_{ir}\in C_{rft}, r\in \sR}{q_{c_{ir}fk}}$) implies zero waiting time as stated in \cref{lemma3} but high charger cost. We use a local search to find a suitable $z_{fk}$ value. In this regard we first calculate the objective value ($\textbf{C}$) for the upper and lower bounds of $z_{fk}$; if the former has a lower objective value, we update $z_{fk}$ by subtracting chargers using \textproc{zlUpdater} (\cref{z_upper_pseudo}). Otherwise (i.e., the latter has a lower objective value), we increase the value of $z_{fk}$ using \textproc{zlUpdater} (\cref{z_lower_pseudo}) until $\textbf{C}$ reaches the minimum. The \textproc{Evaluator} function from (\cref{Eval_pseudo}) updates the value of $z_{fk}$. 

\begin{lemma}\label[lem]{lemma2}
For a given $f\in\sF$ and $k\in \sK{}$, $z_{fk}^*\in [1$, $\sum_{c_{ir}\in C_{rft}, r\in \sR}{q_{c_{ir}fk}}]$, if 
    $\sum_{c_{ir}\in C_{rft}}{q_{c_{ir}fk}}>1$; otherwise, $z_{fk}^* = 0$.
\end{lemma}

\begin{lemma}\label[lem]{lemma3}
For a given $f\in\sF$ and $k\in \sK{}$, $w_{c_{ir}fk} = 0$, if 
    $z_{fk} = \sum_{c_{ir}\in C_{rft}}{q_{c_{ir}fk}}$.
\end{lemma}

In each step of \textproc{Evaluator}, given $q_{c_{ir}fk}$ and $z_{fk}$, we calculate other variables using \textproc{LowerLevelEvaluator} (\cref{LLE_pseudo}). In \textproc{LowerLevelEvaluator}, we first calculate the recharging time ($u_{c_{ir}fk}$) using \textproc{uCalculation} (\cref{u_pseudo}). We assume that trucks recharge at a facility for a duration sufficient to complete the trip if it is less than the maximum battery capacity minus the battery's current level; otherwise, they recharge to full capacity. Next, we calculate the wait time ($w_{c_{ir}fk}$) using \textproc{wCalculation} (\cref{w_pseudo}). To do so, we follow the first-come-first-served rule: that is, vehicles recharge at a facility in the order of their arrival times. Given $q_{c_{ir}fk}$, $u_{c_{ir}fk}$, $w_{c_{ir}fk}$, and $z_{fk}$, calculation of other variables and therefore $\textbf{C}$ is straightforward.

\subsection{A hybrid solution approach supported by the genetic algorithm}\label[sec]{hybrid_sec}

The GA initializes with estimating values for $q_{c_{ir}fk}$ variables. A hybrid solution approach can be formed by feeding these GA-generated $q_{c_{ir}fk}$ solutions into an MILP solver as a constraint set. Therefore, the hybrid approach finds optimal solutions for fixed $q_{c_{ir}fk}$ decisions. This is especially useful because many solutions can be investigated in parallel. In the following section, the performance of this solution approach will be compared with that of the GA.

\section{Case studies}\label[sec]{case_studies} 
This section describes the details of our experimental design and data in \cref{Data structure}, demonstrates the performance and limitations of the GA and the hybrid methods in \cref{computational_exp}, illustrates the impact of the time step duration on system cost in \cref{T_Delta_impact}, and provides key managerial and policy insights along with extensive sensitivity analyses aiming to identify crucial parametric levers in \cref{charger_cost_impact}--\cref{B_bar_cost_impact}.

\subsection{Design of experiments}\label[sec]{Data structure}

We conduct numerical experiments using the Chicago metropolitan area. E-commerce daily demand and road network data are obtained from POLARIS
\citep{auld2016polaris}. We utilize the framework developed by \cite{cokyasar2022time} to form parcel delivery truck routes. For a given set of customer and depot locations and other parameters (e.g., the operational time during a day and vehicle capacity) the framework yields vehicle routes that are the sequences of customers to be visited. Customers in these routes are aggregated at clusters using the RCP solution model in \citet{cokyasar2023}.

\begin{figure}[ht!]
  \begin{center}
    \includegraphics[width=.7\linewidth]{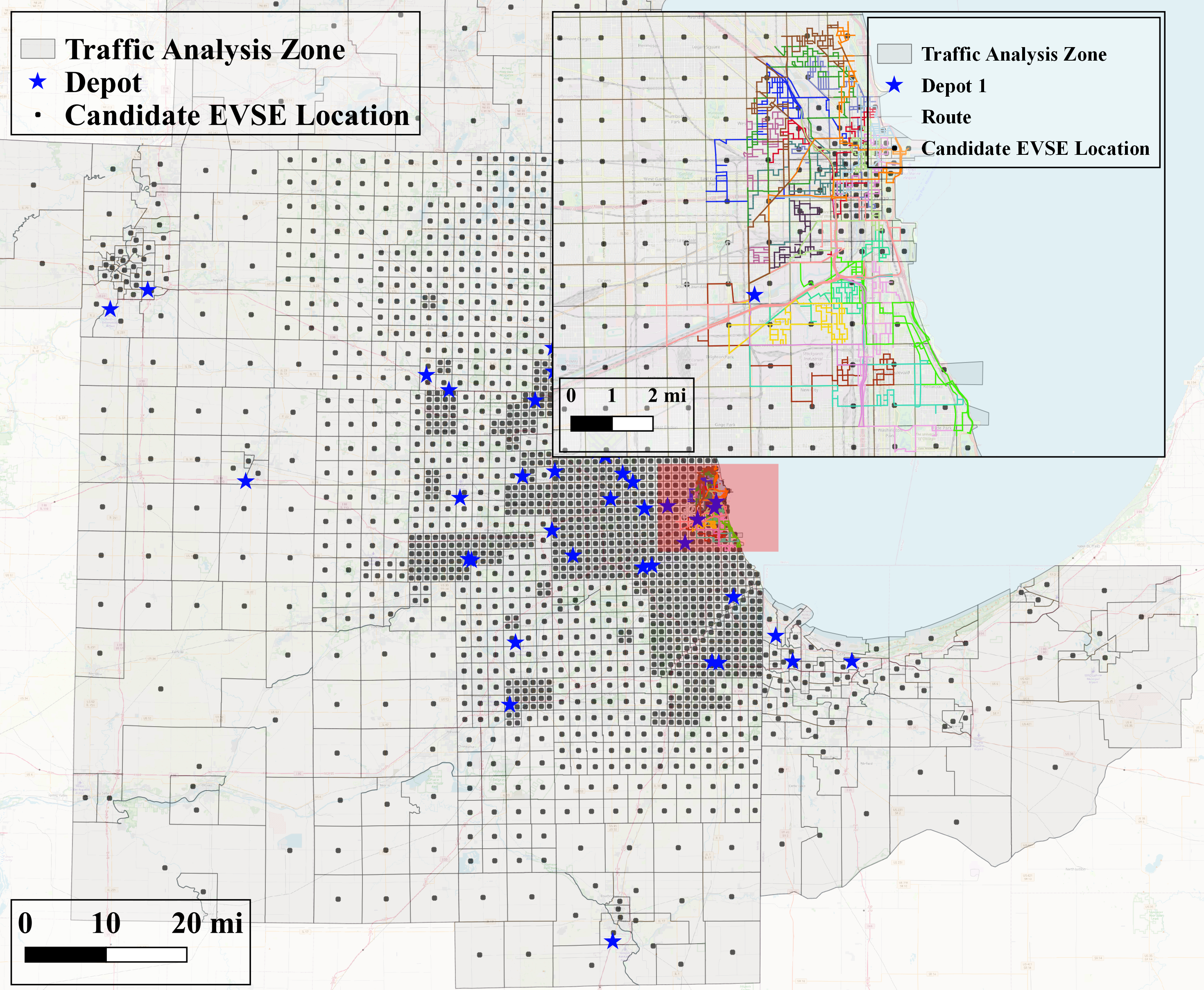}
  \end{center}
  \caption{Illustrative EVSELCA problem instance in the Chicago metropolitan area. The main figure shows traffic analysis zones and depots used for e-commerce delivery in the area as well as candidate charging facility locations. The inset depicts the study region including links of routes and a depot of those routes. Links in the inset are color-coded, and each color indicates a specific route.\label[fig]{route_map}}
\end{figure}

\cref{route_map} depicts an example problem layout. In this figure, the centroids of traffic analysis zones (TAZs), defined by metropolitan area organizations, represent candidate charging facility locations. In the experiments we use 20 routes depicted and their feeder depot, called depot 1. The depot serves 4,580 customers, and customers are aggregated at 89 clusters. Unless otherwise noted, in our experiments we consider four candidate locations, namely, the depot and the three closest centroids to the depot.

The problem parameters used, based on the literature \citep{energy_cost, operation_time, Ellis2017, Davatgari2021, Battery_capacity, Smith2015}, are summarized in \cref{param_value}. Aside from these parameters, we estimate $T_{c_{ir}c_{jr}}^\tau$, $T_{c_{ir}}^\kappa$, and $T_{c_{ir}f}^\delta$ using Manhattan distances and assuming a constant truck speed of 30 mph \citep{speed}. Unless otherwise stated, we let $T^\Delta=15$ minutes. As with most studies in the literature, we consider three charger types with varying powers \citep{Liu2017, 6280677}. The average time required to power up a battery for 100 miles of range gain and charger installation costs are shown in \cref{Charger_config}, which are derived from the literature \citep{Bennett2021}. All costs in the objective function are converted to USD per day. To do so, the lifespans of chargers and facilities are set to 10 and 40 years, respectively \citep{Bennett2021}.

\begin{table}[!htbp]\caption{Parametric values used.}\label[tab]{param_value}
    \footnotesize
    {\begin{tabular*}\textwidth{c@{\extracolsep{\fill}}cccccccc}
        \toprule
        \multicolumn{1}{c}{\begin{tabular}[c]{@{}c@{}}$C_f^\phi$\\ (USD/day)\end{tabular}}
        &\multicolumn{1}{c}{\begin{tabular}[c]{@{}c@{}}$T_{c_{ir}}^\kappa$\\ (minute)\end{tabular}}
        &\multicolumn{1}{c}{\begin{tabular}[c]{@{}c@{}}$\overline{T}$\\ (hour)\end{tabular}}
        &\multicolumn{1}{c}{\begin{tabular}[c]{@{}c@{}}$\overline{B}$\\ (minute)\end{tabular}}
        &\multicolumn{1}{c}{\begin{tabular}[c]{@{}c@{}}$B_r^\iota$\\ (minute)\end{tabular}}
        &\multicolumn{1}{c}{\begin{tabular}[c]{@{}c@{}}$B_r^\omega$\\ (minute)\end{tabular}}
        &\multicolumn{1}{c}{\begin{tabular}[c]{@{}c@{}}$C_k^\xi$\\ (USD/kWh)\end{tabular}}
        &\multicolumn{1}{c}{\begin{tabular}[c]{@{}c@{}}$C^\rho$\\ (USD/mile)\end{tabular}}\\
        \midrule
        35&2&14&200&200&160& 0.43& 1.377\\
        \bottomrule
    \end{tabular*}}
\end{table}

\begin{table}[!htbp]\caption{Charger configuration.}\label[tab]{Charger_config}
\footnotesize
{\begin{tabular*}\linewidth{c@{\extracolsep{\fill}}ccccc}
\toprule
\multicolumn{1}{c}{\begin{tabular}[c]{@{}c@{}}Charger\\ type\end{tabular}}
&\multicolumn{1}{c}{\begin{tabular}[c]{@{}c@{}}Power\\ (kW)\end{tabular}}
&\multicolumn{1}{c}{\begin{tabular}[c]{@{}c@{}}Added driving\\ range (mile)\end{tabular}}
&\multicolumn{1}{c}{\begin{tabular}[c]{@{}c@{}}Added driving\\ charging \\ time (minute) \end{tabular}}
&\multicolumn{1}{c}{\begin{tabular}[c]{@{}c@{}}Cost\\ (USD$\times 1000$)\end{tabular}}\\
\midrule
Basic& 50& 100& 265& 73\\
Moderate& 180& 100& 88& 157\\
Fast& 360& 100& 29& 228\\
\bottomrule
\end{tabular*}}
\end{table}

\subsection{Computational performance of the GA and the hybrid methods}\label[sec]{computational_exp}

We analyze the computational performance of the GA and the hybrid methods and compare them with the MILP solved by Gurobi using $|\sR|$ and $|\sF|$ as problem size determinant levers. The three solution approaches are given the same clusters as an input to make the solutions comparable. A testbed of instances was generated utilizing simulated data from a depot in the Chicago metropolitan area serving e-commerce deliveries. As a baseline, the parametric design provided in \cref{Data structure} was utilized, and we selected three candidate charging facility locations closest to the depot and three random routes out of 20 that the depot serves. While keeping the three charging facility locations and the depot as candidates, a subset of routes $|\sR|\in \{3, 6, 10, 15, 20\}$ were randomly (following a uniform distribution for the selection probability) selected to generate 20 problem instances for each number of routes. Another 20 instances for each of $|\sF|\in \{3, 6, 10, 15, 20\}$ were generated such that the closest charging facility locations to the depot were chosen as candidates, and the three routes in the baseline were used. These 200 instances were then solved by using the three methods with a limit of 600 seconds of computational time per instance. Some of the instances were initially solved without a time limit to observe the impact of the time limit choice. After many compute hours, we did not observe considerable improvement in the solution quality compared to the solutions obtained at the time limit. In the GA, we solved each instance five times and provide statistics of the performances. All computations were carried out on an Intel{\textsuperscript \textregistered} Xeon{\textsuperscript \textregistered} Gold 6138 CPU @2.0 GHz workstation with 128 GB of RAM and 40 cores. Problem instances were solved by using
the Python 3.8.8 interface to the commercial solver Gurobi 10.0~\citep{gurobi}. 

\cref{comp_res} reports the computational performance of the MILP model solved via Gurobi, the GA, and the hybrid solution approaches. The first columns specify the scenario. The MILP columns denote the number of instances that could be solved (i.e., built and reported a feasible solution within 600 seconds) and the number of instances for which optimality was reached, respectively. In the GA, since each instance was solved five times, the maximum number of instances that could be possibly solved was 100 for each scenario. We see that all approaches were unable to produce a feasible solution within the time limit in some runs. GA columns indicate the minimum, maximum, average, and standard deviation of the percent gap between the best objective of the Gurobi-reported solution and the best solution found in the five GA runs. A negative average percent indicates that the GA's best solutions were better than that of the MILP. Hybrid columns follow a similar presentation approach for the hybrid approach.

\begin{table}[!htbp]\caption{\replace{}{Summary of computational performance of the three solution approaches.}}\label[tab]{comp_res}
    \footnotesize
    {\begin{tabular*}\textwidth{cccccccccccccc}
        \toprule
      \multicolumn{2}{c}{Scenario} &
      \multicolumn{2}{c}{MILP} &
      \multicolumn{5}{c}{GA} &
      \multicolumn{5}{c}{Hybrid} \\
      \cmidrule(l{0.5em}r{0.5em}){1-2}
      \cmidrule(l{0.5em}r{0.5em}){3-4} 
      \cmidrule(l{0.5em}r{0.5em}){5-9} 
      \cmidrule(l{0.5em}r{0.5em}){10-14} 
      \multirow{2}{*}[-0.4em]{\begin{tabular}[c]{@{}c@{}}\# \\ $\norm{\sR}$\end{tabular}} &
      \multirow{2}{*}[-0.4em]{\begin{tabular}[c]{@{}c@{}}\# \\ $\norm{\sF}$ \end{tabular}} &
      \multirow{2}{*}[-0.4em]{\begin{tabular}[c]{@{}c@{}}\# \\ Solved\end{tabular}} &
      \multirow{2}{*}[-0.4em]{\begin{tabular}[c]{@{}c@{}}\#\\ Opt\end{tabular}} &
      \multirow{2}{*}[-0.4em]{\begin{tabular}[c]{@{}c@{}}\# \\ Solved\end{tabular}} &
      \multicolumn{4}{c}{$\Delta$ Diff (\%)} &
      \multirow{2}{*}[-0.4em]{\begin{tabular}[c]{@{}c@{}}\# \\ 
      Solved\end{tabular}} &
      \multicolumn{4}{c}{$\Delta$ Diff (\%)} \\
      \cmidrule(l{0.8em}r{0.8em}){6-9}
      \cmidrule(l{0.8em}r{0.8em}){11-14}
    \multicolumn{2}{c}{} &&   &   &
  \multicolumn{1}{l}{Min} &  Max &  Avg &  Std &
   &  Min &
  \multicolumn{1}{l}{Max} &  \multicolumn{1}{l}{Avg} &  Std \\
      \midrule
    \multirow{5}{*}{}3 & 3 & 20 & 19 & 92 & 0    & 19.9   & 2.1    & 5.3  & 17 & 7    & 30   & 19   & 9.7 \\
       6 & 3 & 20 & 0  & 99 & -0.5 & 30   & 6.6  & 7.6 & 12 & 4.1  & 18   & 9.7  & 4.8 \\
       10 & 3 & 16 & 0  & 95 & -6.7 & 6.7  & 0.8  & 3.7 & 3  & -3   & 2.5  & 0.1  & 2.3 \\
       15 & 3 & 6  & 0  & 94 & -10  & 3.8  & -3.5 & 4.5 & 1  & -5.8 & -5.8 & -5.8 & 0   \\
       20 & 3 & 4  & 0  & 89 & -5.2 & -3.2 & -3.8 & 0.8 & 0  & -    & -    & -    & -   \\
    \cmidrule(l{0.5em}r{0.5em}){1-14}
    \multirow{5}{*}{}3 & 3  & 20 & 19 & 88 & 0    & 19.8   & 1.3  & 4.5  & 17 & 8.5  & 29   & 21   & 7.3 \\
       3 & 6  & 20 & 17 & 91 & 0    & 19.8   & 3.2  & 6.1  & 17 & 0    & 31   & 24   & 7.9 \\
       3 & 10 & 20 & 8  & 90 & -0.1 & 25   & 3.9  & 7.8 & 17 & 0    & 27   & 13   & 9.4 \\
       3 & 15 & 20 & 1  & 94 & 0    & 11   & 1.1  & 2.5 & 17 & 0    & 10   & 5.9  & 2.9 \\
       3 & 20 & 20 & 1  & 89 & -1.4 & 15.5 & 1.9  & 3.4 & 17 & 0    & 26   & 5.4  & 5.6\\
        \bottomrule
    \end{tabular*}}
    \replace{}{Note: $\Delta$ Diff is calculated by one minus the division of the best solution obtained by the corresponding method to the best objective reported by the MILP solver.}
\end{table}

In \cref{comp_res}, \replace{}{we acknowledge the inherent limitations regarding the reliability of the $\Delta$ Diff for scenarios where Gurobi failed to converge to optimality (e.g., $|\sR|\in \{6, 10, 15, 20\}$) within the provided time limit. Note that small gap percentages observed should not be interpreted as indicators of good performance, just as large gap percentages observed in the GA and the hybrid approaches should not necessarily be interpreted as indicators of poor performance. These percentages primarily derive from non-optimal solutions obtained from the optimizer. However, we can confidently rely on comparisons to the optimal solutions achieved by Gurobi (namely, $|\sF|\in\{3, 6, 10, 15, 20\}$). It is worth mentioning that} an increase in $|\sR|$ impacts the problem difficulty more than does an increase in $|\sF|$. Using the MILP through a solver can address only small problems. The GA performs better than the hybrid approach. \replace{Note that large gap percentages in the GA and the hybrid should not be a sign of poor performance because these percentages are  based mostly on nonoptimal solutions obtained from the optimizer.}{}

\begin{figure*}[!htbp]
    \centering
    \subfloat[\centering \label{fig:time_res_a}]{%
        \includegraphics*[width=0.495\textwidth,height=\textheight,keepaspectratio]{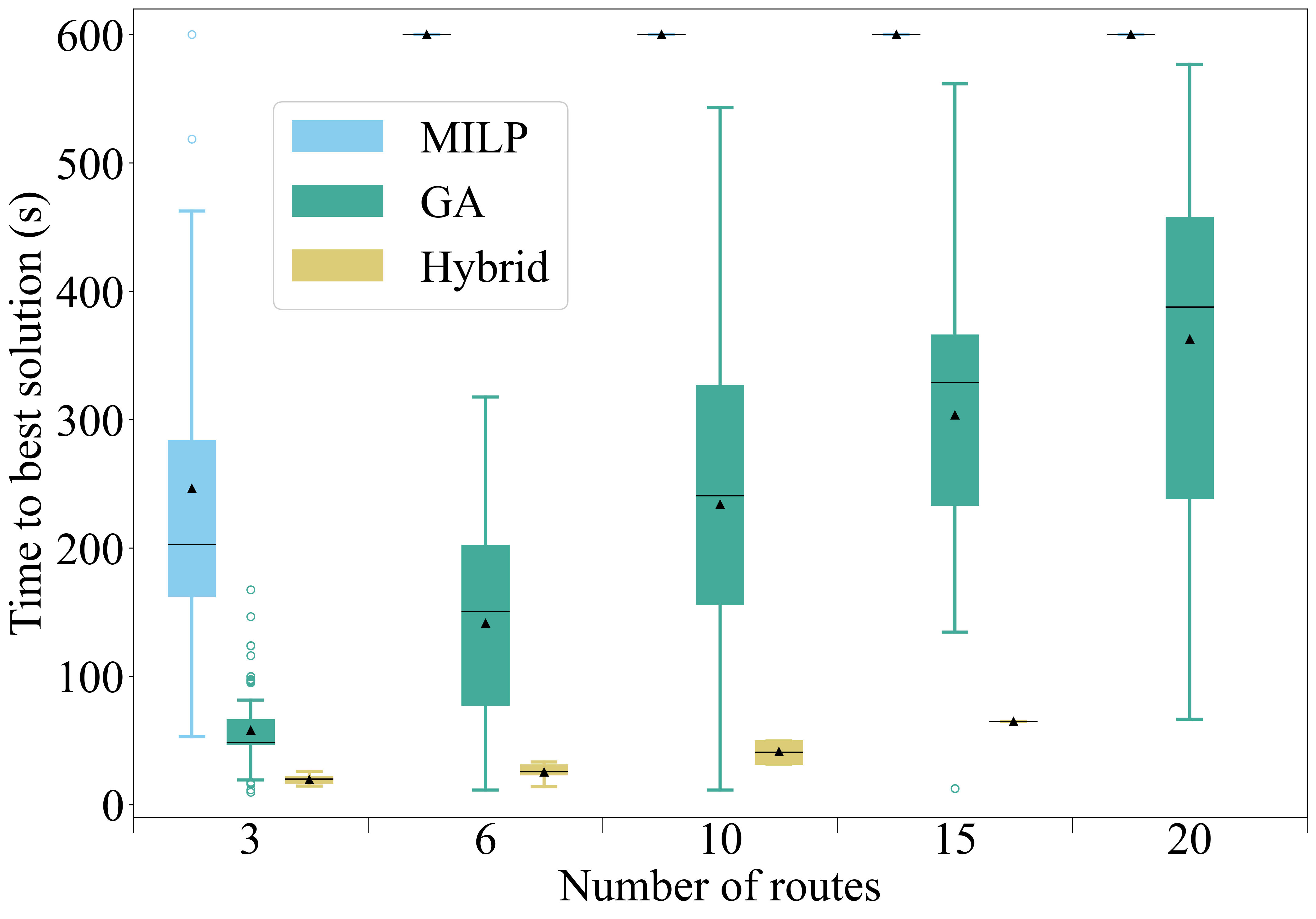}}
    \subfloat[\centering \label{fig:time_res_b}]{\includegraphics*[width=0.495\textwidth,height=\textheight,keepaspectratio]{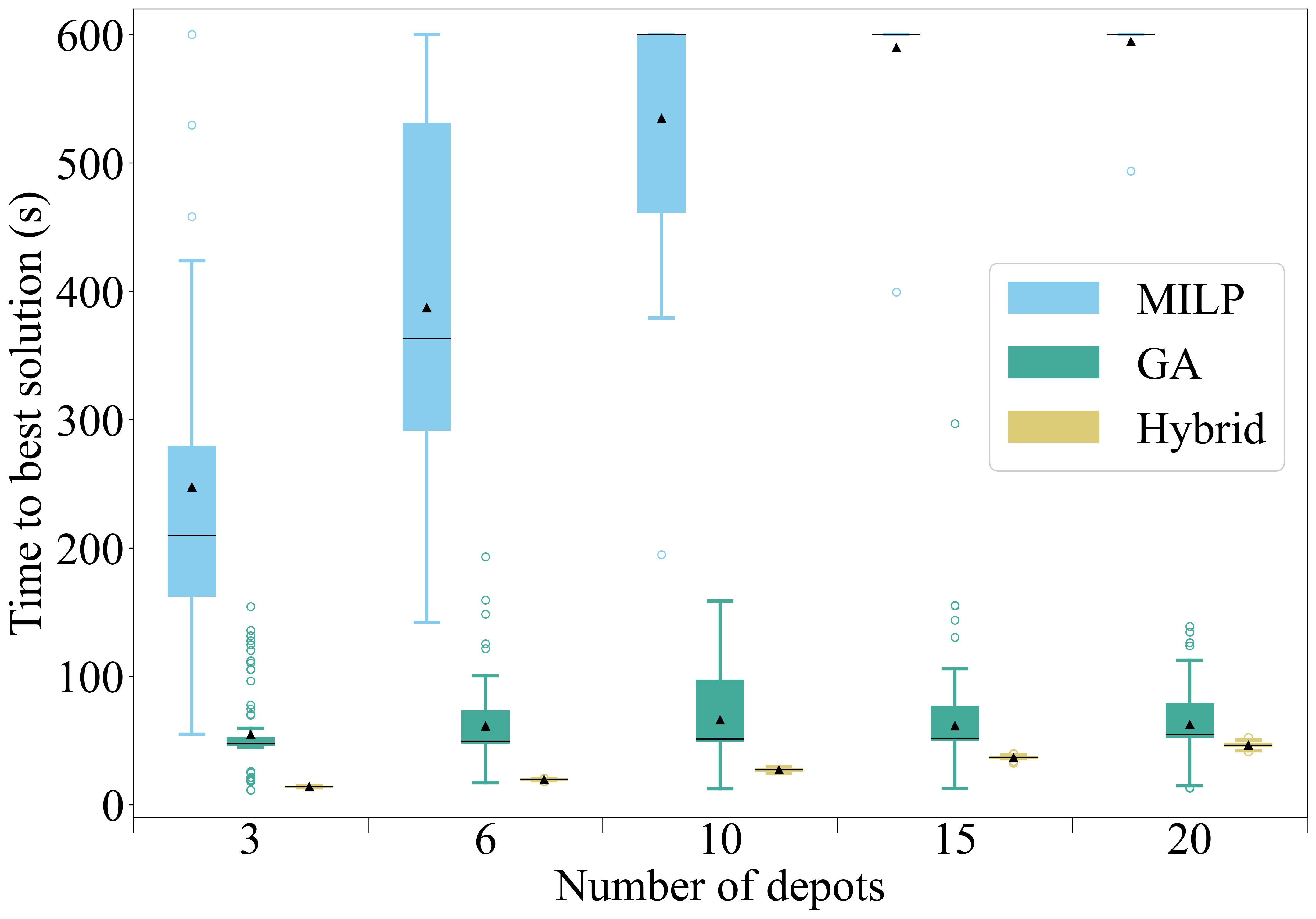}}
    \caption{Time-to-best-solution statistics for the three solution approaches.} \label[fig]{time_res}
\end{figure*}

\cref{time_res} shows the solution time to the best solution per scenario and solution approach. The hybrid approach is the fastest in finding a solution in most scenarios. The time-to-solution comparison between an increasing number of routes and an increasing number of depots supports the claim that $|\sR|$ is a key metric in problem difficulty. In \cref{fig:time_res_a} supported by \cref{comp_res}, we can observe that 6--20 routes scenarios were not solved to optimality within 600 seconds. \cref{fig:time_res_a} and \cref{fig:time_res_b} show that the GA provides quick solutions that are indeed not far off from the MILP (see \cref{comp_res} for a quality comparison).

\subsection{Impact of time step duration on system cost}\label[sec]{T_Delta_impact}

\begin{figure*}[!htbp]
    \centering
    \includegraphics*[width=0.50\textwidth,height=\textheight,keepaspectratio]{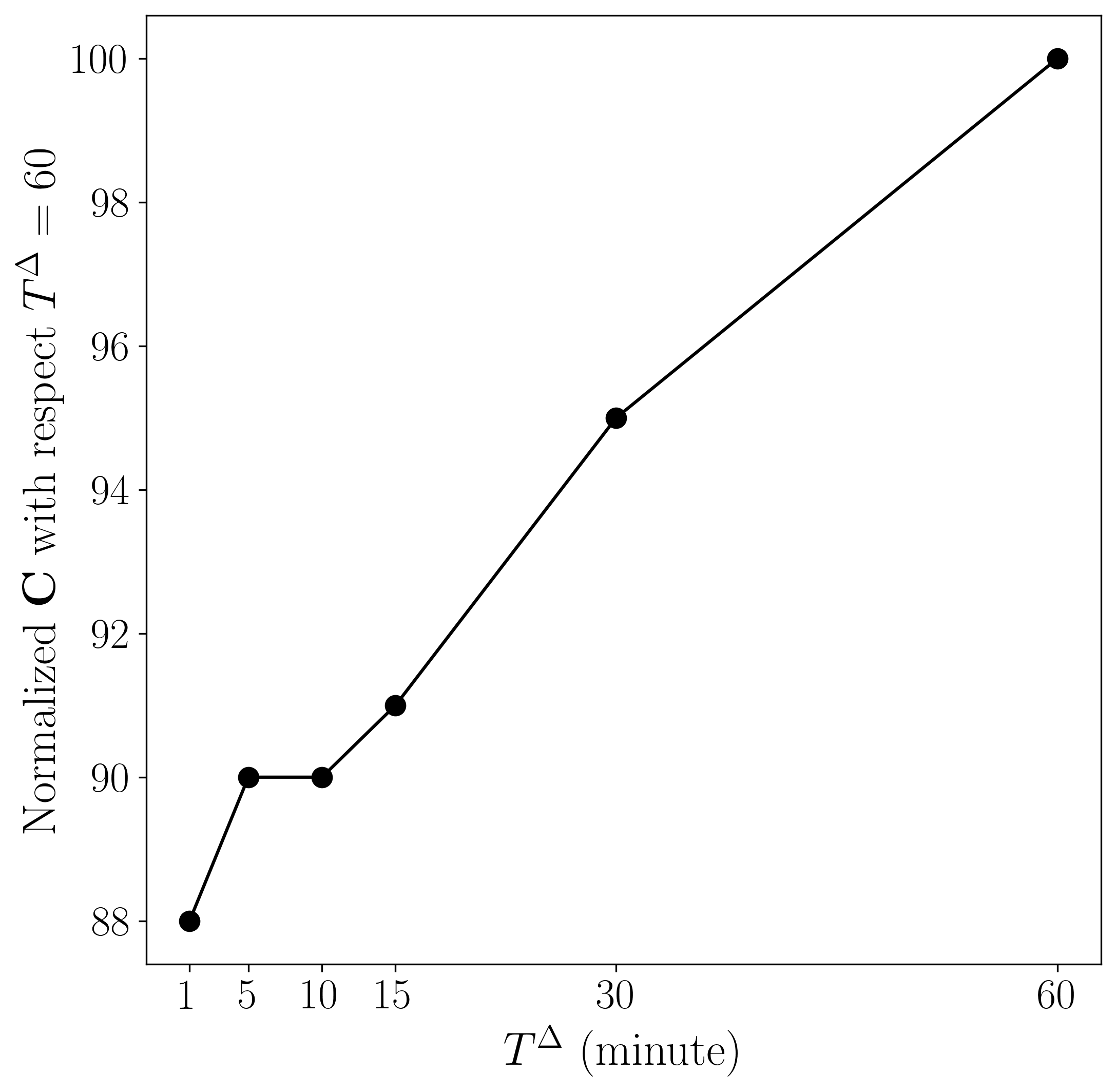}
    \caption{Impact of $T^\Delta$ on $\textbf{C}$.} \label[fig]{T_Delta_fig}
\end{figure*}

The duration of the time step, $T^\Delta$, discretizes the time and control charger availability. Understandably, it can have a considerable impact on the solution quality. It is clear that a large $T^\Delta$ could lead to a higher system cost as more chargers would be needed if an idle charger is shown to be occupied. For instance, $T^\Delta=60$ minutes will label a charger unavailable for an hour even if a vehicle uses the charger for a portion of this time. To assess the impact of $T^\Delta$ on $\textbf{C}$, we consider 1, 5, 10, 15, 30, and 60 minutes as values for the time step duration, and we use the data from depot 1 and its 20 routes. Using the GA method with a time limit of 600 seconds, we solve each instance 20 times and retrieve the minimum $\textbf{C}$. We assume $T^\Delta=60$ as our baseline scenario by setting its $\textbf{C}=100$ and normalize $\textbf{C}$ of other scenarios accordingly. The resulting comparison is demonstrated in \cref{T_Delta_fig}. We can observe that $T^\Delta=1$ yields 12\% lower $\textbf{C}$ compared to that of $T^\Delta=60$.

\subsection{Impact of charger costs on system cost components}\label[sec]{charger_cost_impact}

\begin{figure*}[!htbp]
    \centering
    \subfloat[\centering Percent shares of cost components in $\textbf{C}$.\label{fig:Charger_cost_donut}]{%
        \includegraphics*[width=0.50\textwidth,height=\textheight,keepaspectratio]{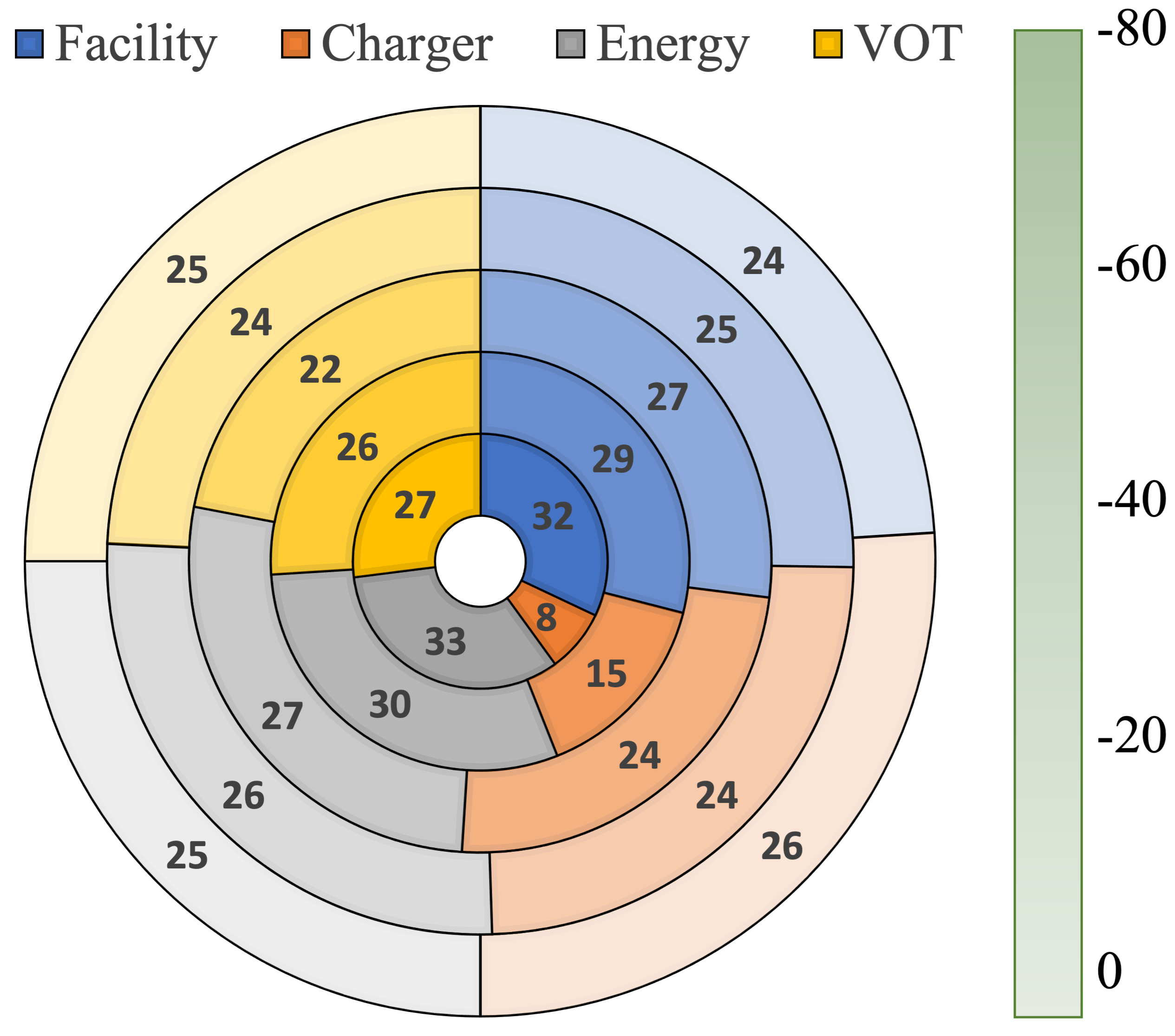}}
    \subfloat[\centering Change in normalized cost and charger allocations.\label{fig:Decrease_in_charger_cost}]{%
        \includegraphics*[width=0.50\textwidth,height=\textheight,keepaspectratio]{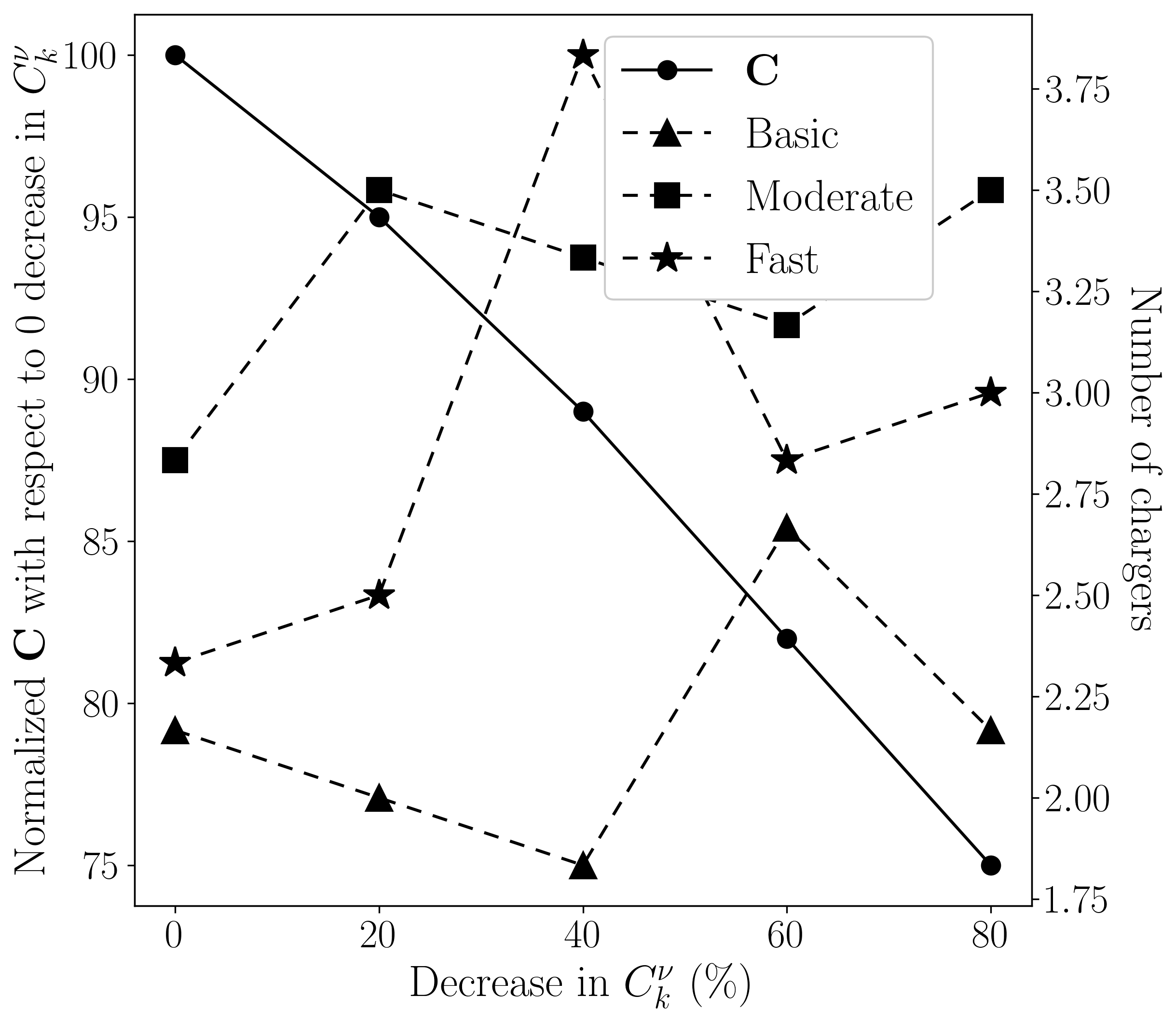}}\\
    \caption{Impact of percent decrease in charger costs.} \label[fig]{charger_cost_fig}
\end{figure*}

Since electrification is a relatively novel technology, charger costs are expected to decrease in the future. Therefore, we analyze the impact of $C_k^\nu$ on the cost components of the objective function \labelcref{obj_fun} and the number of chargers allocated by type. To this end, we define \textit{facility costs} by $\sum_{f\in\sF}C_f^\phi y_f$, \textit{charger costs} by $\sum_{f\in\sF, k\in\sK}C_k^\nu z_{fk}$, \textit{energy costs} by $\sum_{c_{ir}\in\sC_r, r\in\sR, f\in\sF_{c_{ir}}, k\in\sK} C_k^\xi u_{c_{ir}fk}$, and \textit{VOT costs} by $\sum_{c_{ir}\in\sC_r, r\in\sR, f\in\sF_{c_{ir}}, k\in\sK} C^\rho (T_{c_{ir}fk} q_{c_{ir}fk}$ $+w_{c_{ir}fk}+u_{c_{ir}fk})$. Five percentages of decrease in $C_k^\nu$ are considered: 0\%, 20\%, 40\%, 60\%, and 80\%. We consider depot 1 and the 20 routes it serves by sampling $|\sR|\in \{10, 12, 14, 16, 18, 20\}$ and $|\sF|\in \{100\}$. We use the GA method, limit the solution time per instance by 600 seconds, and solve each instance 20 times to obtain the best solution. We conducted 600 ($6\times 5\times 20$) runs for this analysis, and we report the statistics of the best solutions out of 20 GA runs by aggregating over $|\sR|$ in \cref{charger_cost_fig}. 

\cref{fig:Charger_cost_donut} shows that the contribution of charger costs into \labelcref{obj_fun} drops from 26\% to 8\% parallel to $C_k^\nu$, while facility and energy costs substantially increase by 8\% (from 24 to 32). Note that -80, -60, -40, and -20 indicate a 80\%, 60\%, 40\%, and 20\% decrease in charger costs, while 0 shows no change. This tiered representation is adopted in the other donut charts as well, and a negative value indicates a decrease, while a positive one refers to an increase in the corresponding parameter.

\cref{fig:Decrease_in_charger_cost} shows the decrease percent in $C_k^\nu$ versus the normalized cost. The normalized cost assumes $\textbf{C}=100$ when the decrease in $C_k^\nu$ is 0 and is calculated accordingly for other instances. In this figure we observe that a large decrease in $C_k^\nu$ increases the number of moderate and fast chargers, although there are fluctuations. The spikes can be a result of  finding nonoptimal solutions through the GA and finding a better solution  by swapping types of chargers along with a drop in $C_k^\nu$. An example of the latter can be observed by seeing the number of chargers for basic and fast moving from 40\% to 60\%: that is, fewer fast chargers are equipped, while more basic chargers are utilized at 60\% compared with 40\%. 

We note that an 80\% decrease in $C_k^\nu$ results in a 25\% drop in $\textbf{C}$. This percentage will indicate the importance of $C_k^\nu$ compared with others analyzed in the following sections.

\subsection{Impact of energy costs on system cost components}\label[sec]{energy_cost_impact}

\begin{figure*}[!htbp]
    \centering
    \subfloat[\centering Percent shares of cost components in $\textbf{C}$.\label{fig:Energy_cost_donut}]{%
        \includegraphics*[width=0.50\textwidth,height=\textheight,keepaspectratio]{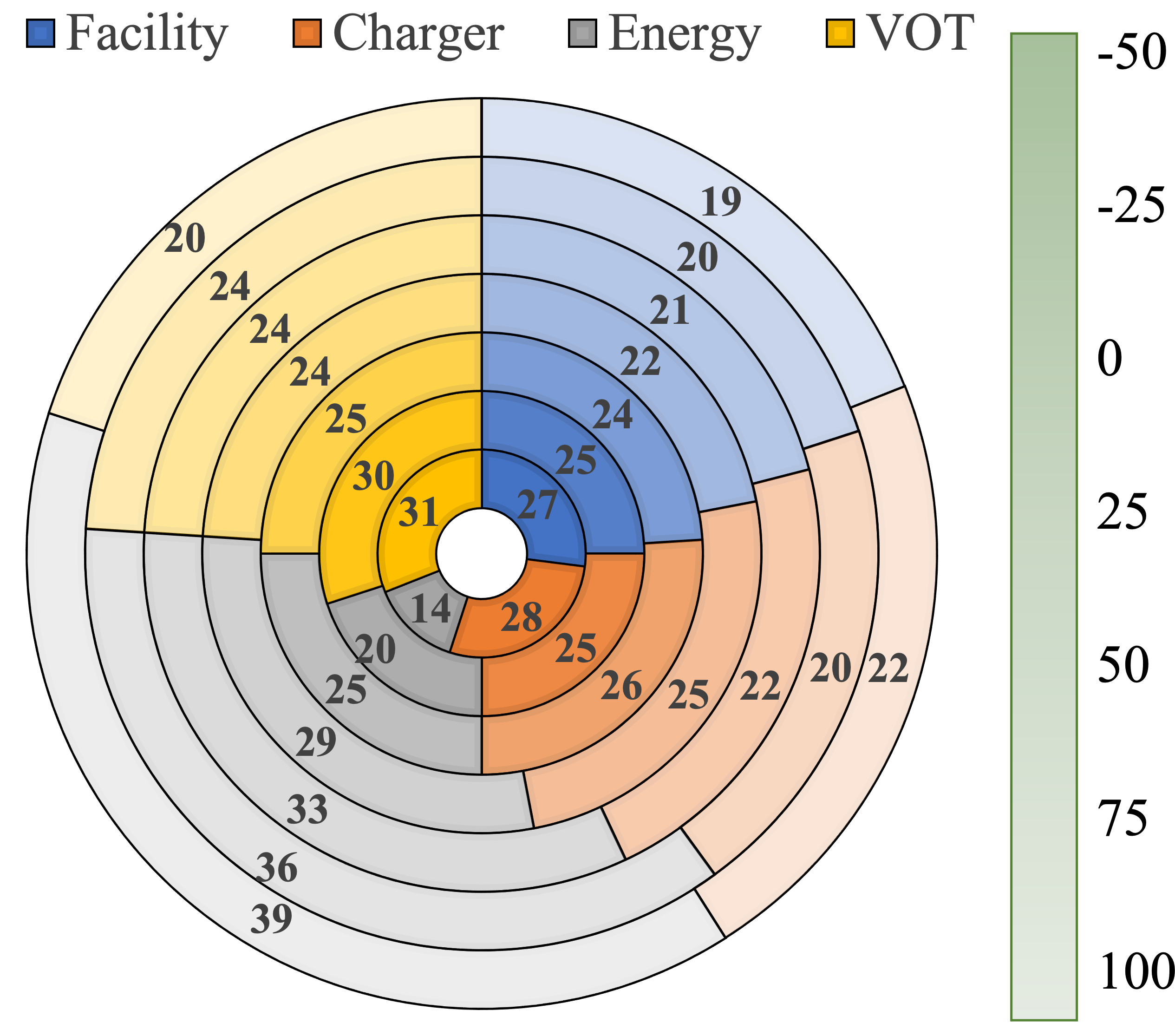}}
    \subfloat[\centering Change in normalized cost and charger allocations.\label{fig:Change_in_energy_cost}]{%
        \includegraphics*[width=0.50\textwidth,height=\textheight,keepaspectratio]{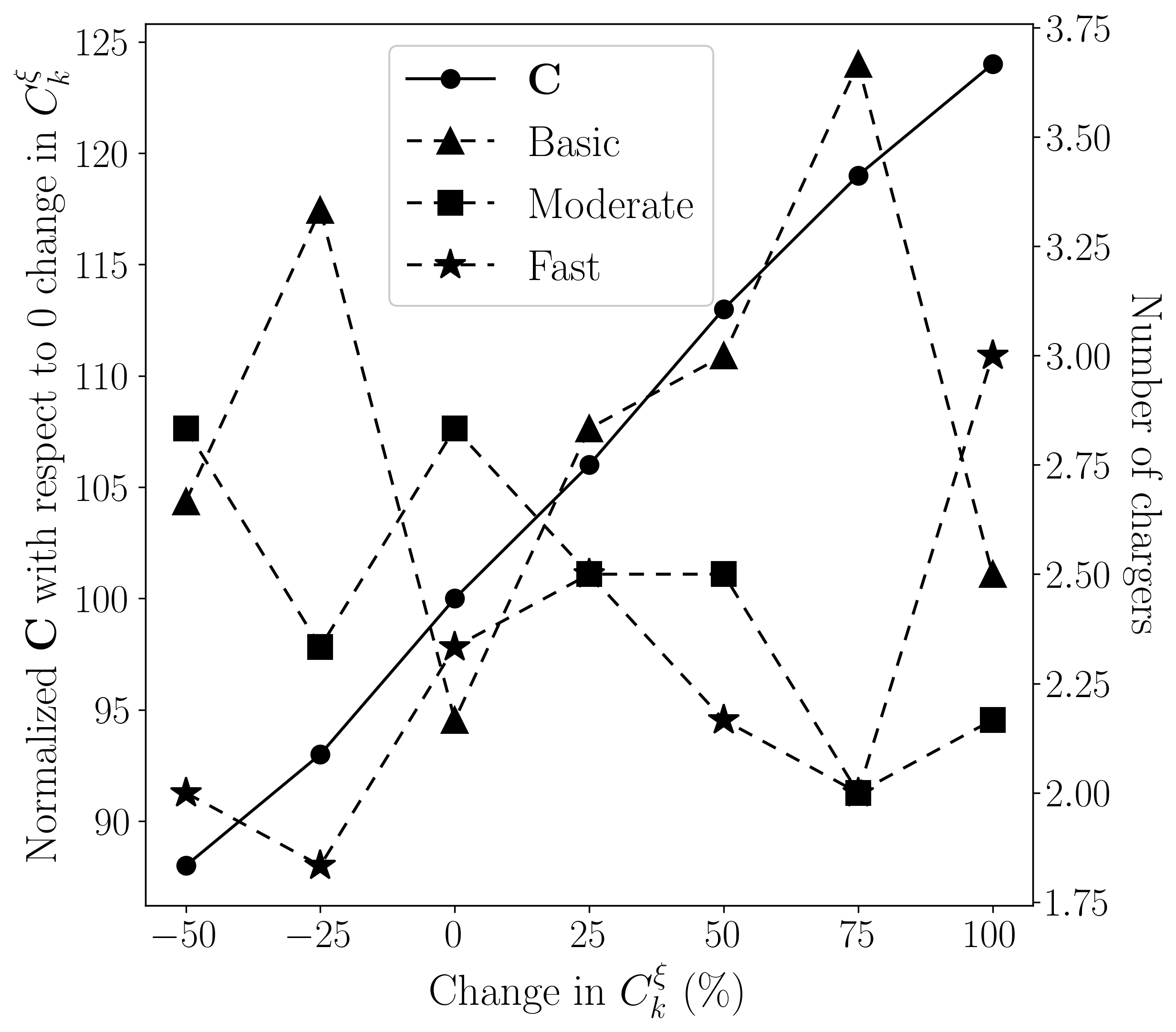}}\\
    \caption{Impact of percent change in energy costs.} \label[fig]{energy_cost_fig}
\end{figure*}

It is not certain how the wider adoption of EVs will impact energy prices. To analyze the impact of energy costs, we consider a percent change of -50, -25, 0, 25, 50, 75, and 100 in $C_k^\xi$. The same settings as in the previous analysis are followed here and will be used in the upcoming sections. Therefore, 840 ($6\times 7\times 20$) runs were conducted, and \cref{energy_cost_fig} provides the summary statistics.

\cref{fig:Energy_cost_donut} shows a substantial jump in VOT costs along with an increase in $C_k^\xi$. In \cref{fig:Change_in_energy_cost}, moving from 0 to 100\% increase in $C_k^\xi$, the number of fast chargers rises, while the numbers of other charger types reduce. 

\subsection{Impact of time value costs on system cost components}\label[sec]{VOT_cost_impact}

\begin{figure*}[!htbp]
    \centering
    \subfloat[\centering Percent shares of cost components in $\textbf{C}$.\label{fig:VOT_cost_donut}]{%
        \includegraphics*[width=0.50\textwidth,height=\textheight,keepaspectratio]{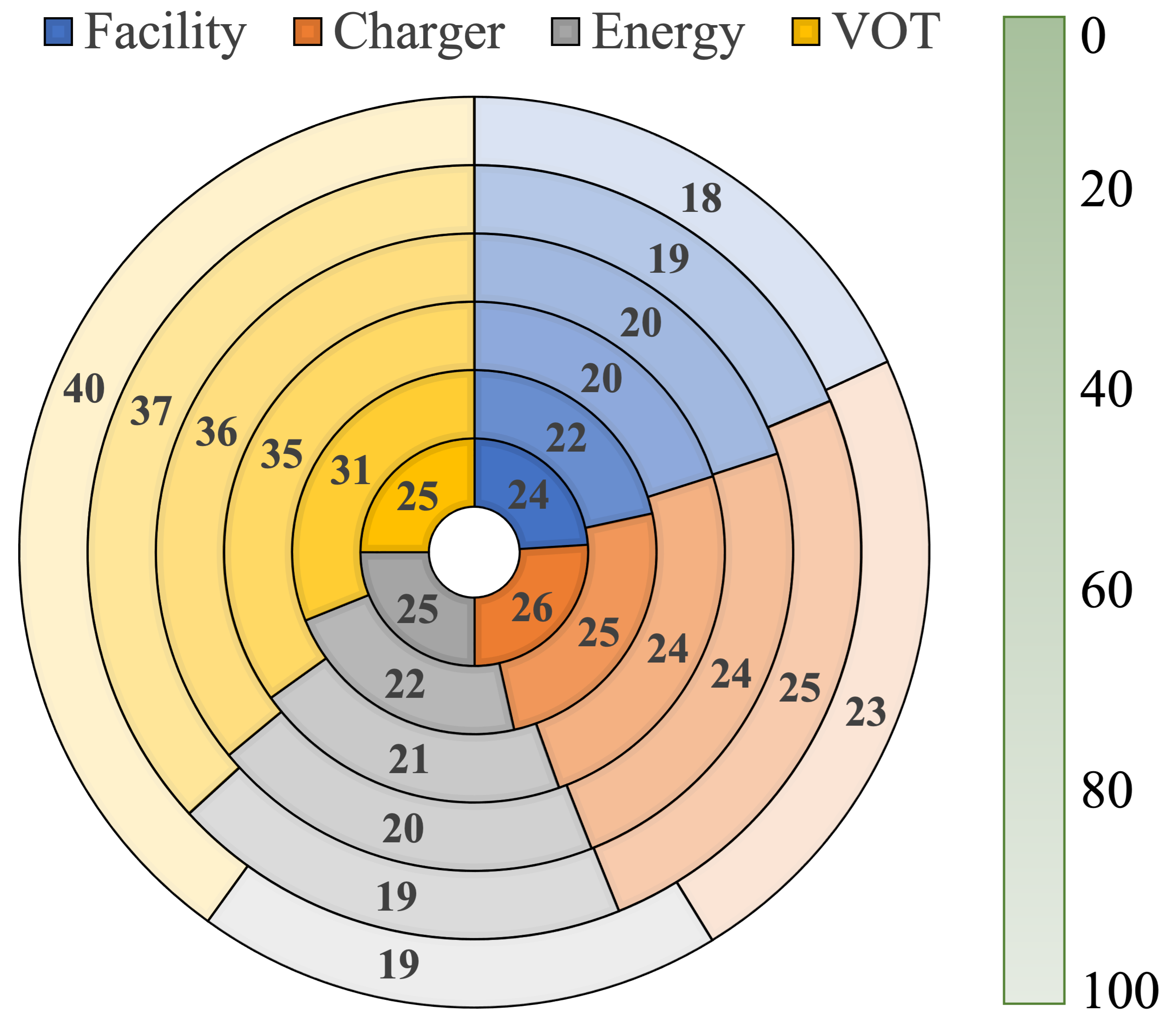}}
    \subfloat[\centering Change in normalized cost and charger allocations.\label{fig:Increase_in_VOT_cost}]{%
        \includegraphics*[width=0.50\textwidth,height=\textheight,keepaspectratio]{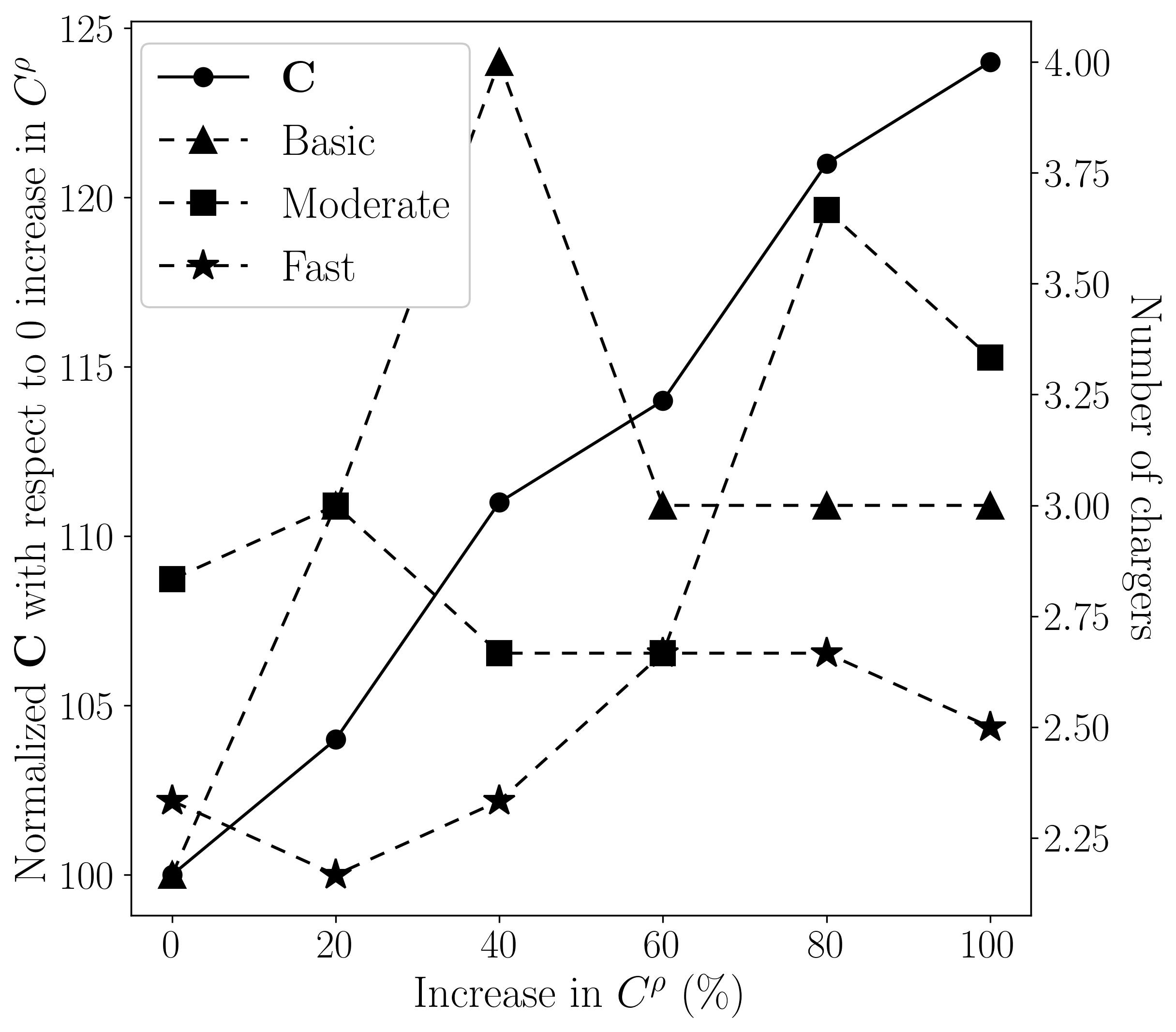}}\\
    \caption{Impact of percent increase in VOT costs.} \label[fig]{VOT_cost_fig}
\end{figure*}

In some industries, VOT may be more important than others. Consumers may be willing to pay more for faster delivery. To analyze how a percent increase of 0, 20, 40, 60, 80, and 100 in $C^\rho$ impact the decisions, we conducted 720 ($6\times 6\times 20$) runs and report summary statistics in \cref{VOT_cost_fig}.

\cref{fig:VOT_cost_donut} illustrates how VOT costs can become  dominant (by 40\%) in \labelcref{obj_fun} when $C^\rho$ is doubled. From \cref{fig:Increase_in_VOT_cost}, we observe that the number of moderate and fast chargers  increases as $C^\rho$ doubles.

\subsection{Impact of battery capacity on system cost components}\label[sec]{B_bar_cost_impact}

\begin{figure*}[!htbp]
    \centering
    \subfloat[\centering Percent shares of cost components in $\textbf{C}$.\label{fig:B_bar_donut}]{%
        \includegraphics*[width=0.50\textwidth,height=\textheight,keepaspectratio]{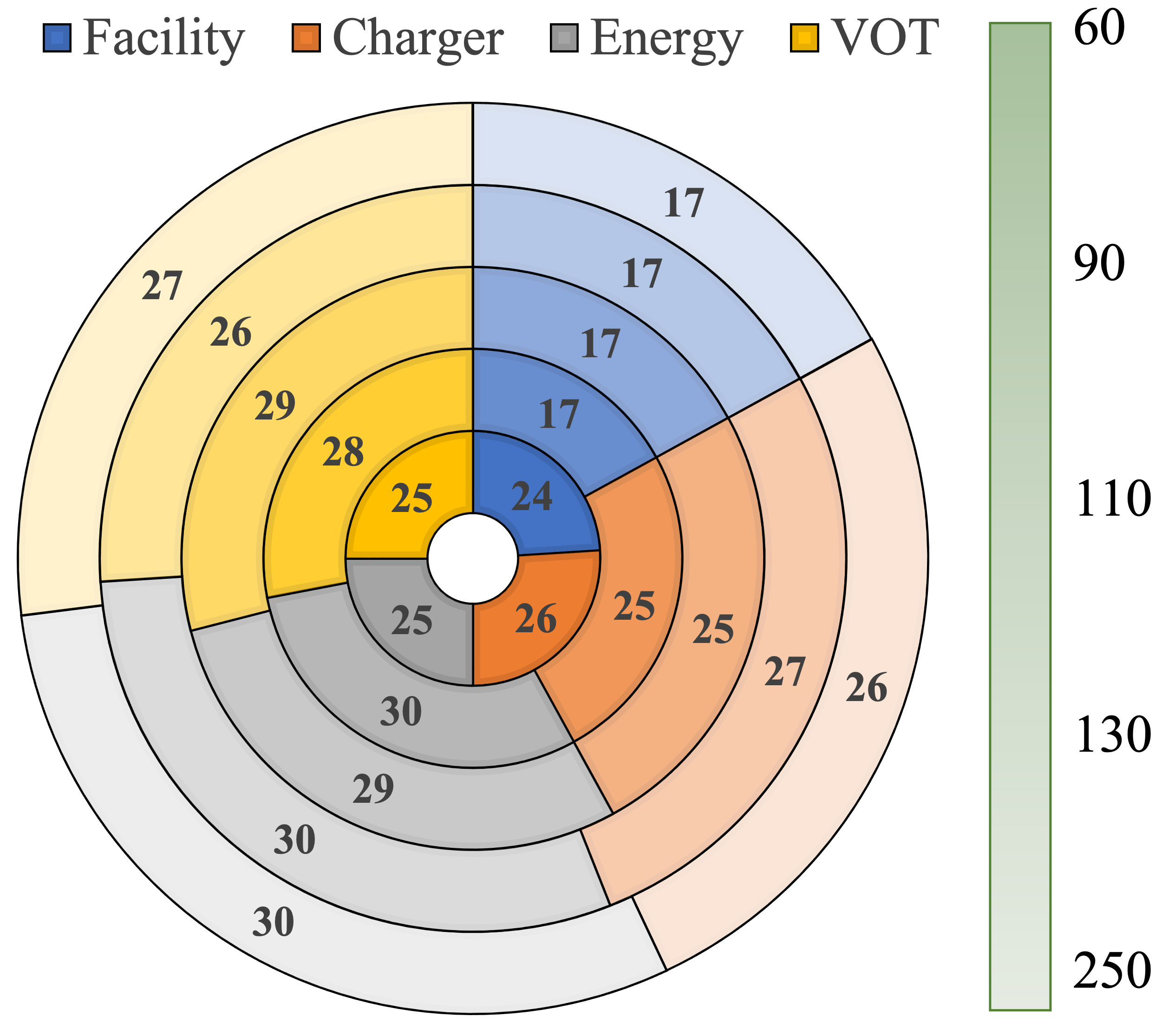}}
    \subfloat[\centering Change in normalized cost and charger allocations.\label{fig:Increase_in_B_bar}]{%
        \includegraphics*[width=0.50\textwidth,height=\textheight,keepaspectratio]{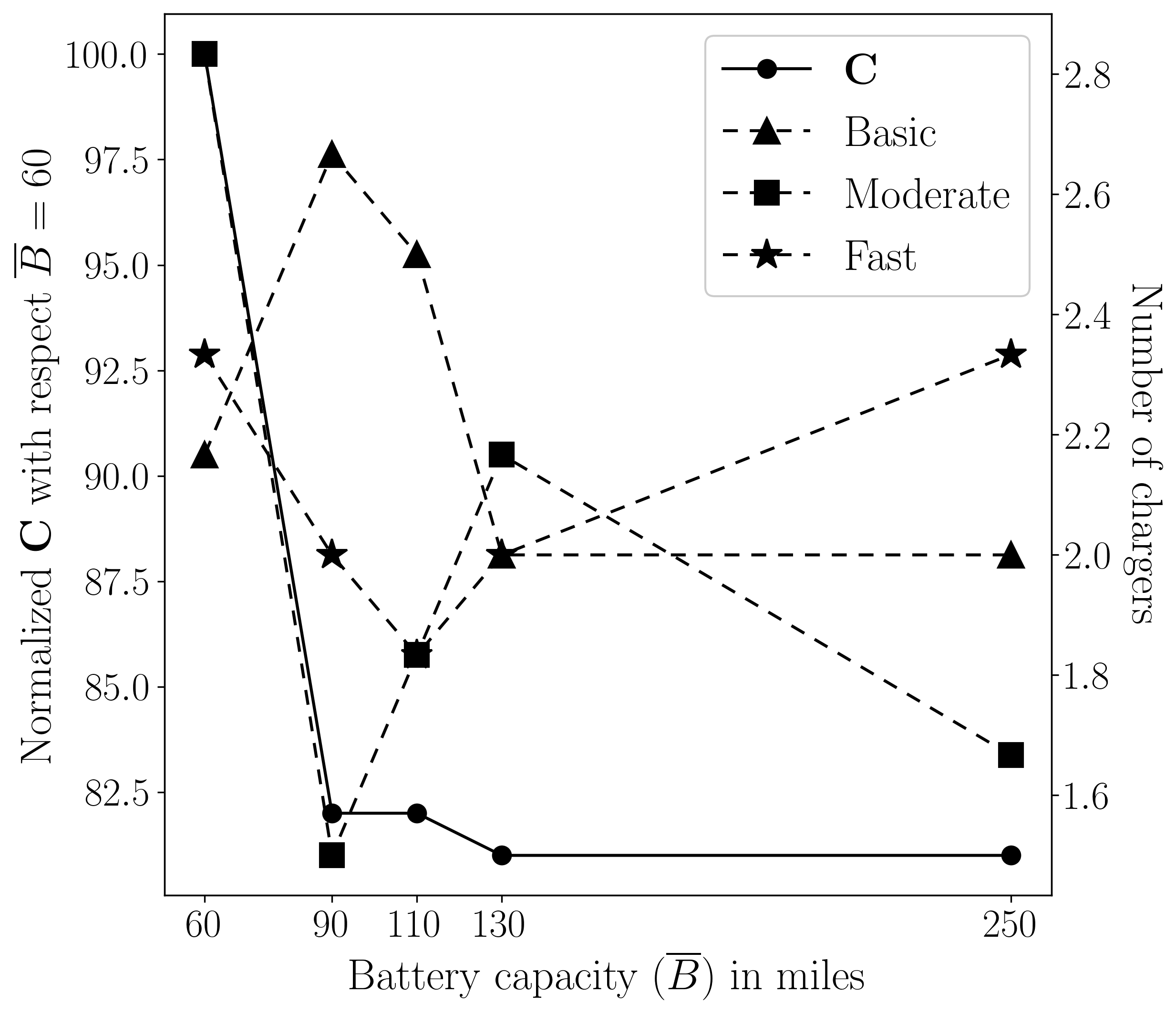}}\\
    \caption{Impact of battery capacity (in miles range).} \label[fig]{B_bar_cost_fig}
\end{figure*}

EV technology is continuously improving, and advancements in battery technology enable longer vehicle ranges. In all analyses, we considered the low end of a 60--130 EV range reported \citep{Battery_capacity} to be conservative. We now analyze EV ranges of 60, 90, 110, 130, and 250 miles. We conducted 600 ($6\times 5\times 20$) runs for this analysis, and the findings are illustrated in \cref{B_bar_cost_fig}.

Increasing the EV range reduces facility costs by enabling vehicles to recharge at more central locations, as shown in \cref{fig:B_bar_donut}. \cref{fig:Increase_in_B_bar} demonstrates that a longer EV range decreases $\textbf{C}$ up to a point. A similar finding was observed in a previous study \citep{cokyasar2022time}.

\section{Conclusion}\label[sec]{conclusion}
This research addresses some challenges in electrifying trucks, particularly in designing the necessary charging infrastructure. Indeed, infrastructure planning for large-scale electrification projects can be complex because of factors including the placement of charging stations and the scheduling of recharging activities. While previous studies have focused on solving the EVLRP to address this issue, the freight industry has a different approach, where they make electrification decisions based on the characteristics of existing routes. However, by prioritizing routes shorter than the EV range and installing chargers only at depots, typically one charger for each EV, they forego the potential for cost savings that can be achieved by optimizing the charging activities. Exploiting and studying the limits of the associated economic opportunities is crucial, given the long-term impacts of strategic location-allocation decisions on short-term routing activities. Our paper seeks to address this gap by providing decision-making models that combine the best of existing studies while also respecting the freight industry's philosophy of electrifying existing routes.

To that end, we formally describe the EVSELCA problem and construct an MILP model that focuses on fixed routes. The MILP model, solved through commercially available solvers that often use branch and bound as a solution method, faces scalability issues, making it impractical for larger-scale problems. To overcome this challenge, we propose a clustering approach that simplifies the problem by grouping customers into clusters and allowing recharging only upon completion of service at these clusters. Clustering is shown to partially address computational difficulties for small-scale problems, yet it is not adequate to fully address the issue. For this reason we develop a metaheuristic solution method using a GA. This approach generates near-optimal solutions within a reasonable time frame, making it possible to apply the model at a large scale. Combining the GA and MILP solvers, we introduce a hybrid solution approach. 

We compare the computational performance of the GA and hybrid methods with the MILP model solved through Gurobi, using the number of routes and charging facility locations as problem-size-determinant levers. The results indicate that the GA outperforms the hybrid method in terms of solution quality, but the hybrid method is faster in finding solutions in most scenarios. The MILP model is suitable for small-scale problems. Overall, the GA provides quick solutions that are close in quality to the optimal solution. Moreover, the findings show that an increase in the number of routes has a greater impact on problem difficulty than does an increase in the number of charging facility locations. \replace{}{Other evolutionary algorithms or nature-inspired metaheuristics could be tested against the GA's performance, but a systematic comparison is out of scope for the present manuscript.}

We investigate the impact of four factors on the EVSELCA problem through a sensitivity analysis. These are charger costs ($C_k^\nu$), energy costs ($C_k^\xi$), VOT ($C^\rho$), and battery capacity ($\overline{B}$). Our key findings are summarized as follows. An 80\% decrease in $C_k^\nu$ results in a 25\% cost reduction. A substantial reduction in VOT cost shares is observed as $C_k^\xi$ increases. The number of moderate and fast chargers increases as $C^\rho$ doubles. Longer EV ranges are beneficial in decreasing the overall cost up to a certain point. Beyond this threshold, longer EV ranges result in only a negligible decrease in the total cost. Our findings indicate that the objective function is the most sensitive to charger costs compared with other factors, while energy and VOT costs are less vital.

The proposed model is intricate and addresses critical concerns of determining the optimal time for EVs to visit charging facilities, selecting suitable facilities, allocating appropriate charging infrastructure, scheduling recharging activities to minimize wait times, and satisfying operational constraints. The strength of the MILP developed in this study is that it addresses all of these concerns. However, the model's dependency on the candidate locations for charging facility placement presents a challenge: that is, changing one candidate location may substantially alter the solutions and their interpretation. To address this, we plan to develop a tool that can mimic the MILP model to quickly find near-optimal solutions for any given set of candidate locations, reducing the time required to solve the problem.

\section*{Acknowledgments}
This material is based upon work supported by the U.S. Department of Energy, Office of Science, under contract number DE-AC02-06CH11357. 
This report and the work described were sponsored by the U.S. Department of Energy (DOE) Vehicle Technologies Office (VTO) under the Systems and Modeling for Accelerated Research in Transportation (SMART) Mobility Laboratory Consortium, an initiative of the Energy Efficient Mobility Systems (EEMS) Program. 
Erin Boyd, a DOE Office of Energy Efficiency and Renewable Energy (EERE) manager, played an important role in establishing the project concept, advancing implementation, and providing guidance. The authors  also thank Aymeric Rousseau, Joshua Auld, Vincent Freyermuth, Hyunseop Uhm, and Olcay Sahin for their continuous support.

\section*{Author contributions}
The authors confirm contribution to the paper as follows: Study conception and
design: A. Davatgari, T. Cokyasar, and A. Subramanyam; data collection: A. Davatgari; analysis and interpretation of results: A. Davatgari and T. Cokyasar; draft manuscript
preparation: A. Davatgari, T. Cokyasar and J. Larson. All authors reviewed the results and
approved the final version of the manuscript.

\appendix

\renewcommand{\thesection}{Appendix}
\section{Proofs}\label{proofs}
\renewcommand{\thetable}{A.\arabic{table}}
\renewcommand{\thefigure}{A.\arabic{figure}}
\renewcommand{\thealgocf}{A.\arabic{algocf}}
\setcounter{table}{0}
\setcounter{figure}{0}
\setcounter{algocf}{0}

\begin{proof}[Proof of Lemma 1] 
  We divide the proof of \cref{lemma2} into two parts: first we prove $z_{fk}^*
  \leq \sum_{c_{ir}\in C_{rft}, r\in \sR}{q_{c_{ir}fk}}$, then $z_{fk}^* \geq
  1$.
  To prove $z_{fk}^* \leq \sum_{c_{ir}\in C_{rft}, r\in \sR}{q_{c_{ir}fk}}$, we
  use proof by contradiction.
  We suppose $z_{fk}^* > \sum_{c_{ir}\in C_{rft}, r\in
  \sR}{q_{c_{ir}fk}}$. 
  From \labelcref{detour_if_decided_to_charge2} we have
  $x_{c_{ir}fkt} \leq q_{c_{ir}fk}$, which if we aggregate both sides of
  inequality on $c_{ir}$, $\sum_{c_{ir}\in C_{rft}, r\in \sR}{x_{c_{ir}fkt}}
  \leq \sum_{c_{ir}\in C_{rft}, r\in \sR}{q_{c_{ir}fk}}$ holds. 
  Then, we can
  conclude that $z_{fk}^* >$ $\sum_{c_{ir}\in C_{rft}, r\in
  \sR}{x_{c_{ir}fkt}}$, which indicates that for a given $f\in\sF$ and
  $k\in\sK$, at all $t\in\sT_{c_{ir}f}$ there exist unused chargers.
  This is in
  contradiction with the optimality of the $z_{fk}^*$, and therefore our
  assumption that $z_{fk}^* > \sum_{c_{ir}\in C_{rft}, r\in \sR}{q_{c_{ir}fk}}$
  is false, so $z_{fk}^* \leq \sum_{c_{ir}\in C_{rft}, r\in \sR}{q_{c_{ir}fk}}$
  holds. 
  To prove $z_{fk}^* \geq 1$, from the condition of the \cref{lemma2},
  we have $\sum_{c_{ir}\in C_{rft}, r\in \sR}{q_{c_{ir}fk}} \geq 1$ and from
  \labelcref{detour_if_decided_to_charge2} we know $q_{c_{ir}fk} \leq
  \sum_{t\in\sT_{c_{ir}f}}x_{c_{ir}fkt}$, which if we aggregate both sides of
  inequality on $c_{ir}$, $\sum_{c_{ir}\in C_{rft}, r\in
  \sR}{\sum_{t\in\sT_{c_{ir}f}}x_{c_{ir}fkt}} \geq $ $\sum_{c_{ir}\in C_{rft},
  r\in \sR}{q_{c_{ir}fk}}$  holds. 
  Then, we can conclude that $\sum_{c_{ir}\in
  C_{rft}, r\in \sR}{\sum_{t\in\sT_{c_{ir}f}}x_{c_{ir}fkt}} \geq 1$.
  This indicates at least for one $t\in\sT_{c_{ir}f}$, $\sum_{c_{ir}\in
  C_{rft}, r\in \sR}x_{c_{ir}fkt} \geq 1$. 
  So from \labelcref{charger_availability}, we
  can conclude that for a given $f\in\sF$ and $k\in\sK$, if $\sum_{c_{ir}\in
  C_{rft}, r\in \sR}{q_{c_{ir}fk}} \geq 1$, then $z_{fk}^* \geq 1$; otherwise,
  $\sum_{c_{ir}\in C_{rft}, r\in \sR}{q_{c_{ir}fk}} = 0$, and from first part
  of the proof we know $z_{fk}^* \leq \sum_{c_{ir}\in C_{rft}, r\in
  \sR}{q_{c_{ir}fk}}$, so $z_{fk}^* = 0$. 
\end{proof}

\noindent\begin{proof}[Proof of Lemma 2] 
  To prove \cref{lemma3}, we aggregate
  \labelcref{detour_if_decided_to_charge2} on $c_{ir}$, and we get
  \[
  \sum_{\substack{c_{ir}\in C_{rft}\\ r\in \sR}}{x_{c_{ir}fkt}} = \sum_{\substack{c_{ir}\in C_{rft}\\ r\in \sR}}{x_{c_{ir}fkt}^{\beta}} - \sum_{\substack{c_{ir}\in C_{rft}\\ r\in \sR}}{x_{c_{ir}fkt}^{\alpha}}.
  \]
  If we suppose that all routes recharge at
  charger $k$ in facility $f$ at time $t$, then we have 
  \[
  \sum_{\substack{c_{ir}\in C_{rft}\\ r\in \sR}}{x_{c_{ir}fkt}^{\beta}} - \sum_{\substack{c_{ir}\in C_{rft}\\ r\in \sR}}{x_{c_{ir}fkt}^{\alpha}} = \sum_{\substack{c_{ir}\in C_{rft}\\ r\in \sR}}{q_{c_{ir}fk}}.
  \]
  Based on the definition of $x_{c_{ir}fkt}^{\alpha}$ and
  $x_{c_{ir}fkt}^{\beta}$, if a route recharges at charger $k$ in facility $f$
  at time $t$ after serving customer $c_{ir}$, $x_{c_{ir}fkt}^{\alpha} = 0$ and
  $x_{c_{ir}fkt}^{\beta} = 1$.
  Therefore, in the case that all routes are
  recharging at the same charger and time, $\sum_{\substack{c_{ir}\in C_{rft}\\ r\in \sR}}{x_{c_{ir}fkt}^{\alpha}} = 0$ and $\sum_{\substack{c_{ir}\in C_{rft}\\ r\in \sR}}{x_{c_{ir}fkt}^{\beta}} = 1$. 

  If we aggregate
  \labelcref{x_alpha_formulation1} and \labelcref{x_alpha_formulation2} on
  $c_{ir}$, we get 
  \begin{align*}
    \sum_{\substack{c_{ir}\in C_{rft}\\ r\in \sR}}{d_{c_{ir}}} + \sum_{\substack{c_{ir}\in C_{rft}\\ r\in \sR}}{T_{c_{ir}f}^\tau {q_{c_{ir}fk}}} + \sum_{\substack{c_{ir}\in C_{rft}\\ r\in \sR}}{w_{c_{ir}fk}}  \leq
    \sum_{\substack{c_{ir}\in C_{rft}\\ r\in \sR}}{t} + \sum_{\substack{c_{ir}\in C_{rft}\\ r\in \sR}}{T^\Delta}-\sum_{\substack{c_{ir}\in C_{rft}\\ r\in \sR}}{\epsilon} \\+ \bigM \sum_{\substack{c_{ir}\in C_{rft}\\ r\in \sR}}\left({1-{q_{c_{ir}fk}}+x_{{c_{ir}fkt}}^{\alpha}}\right)
  \end{align*}
  and
  \begin{align*}
    \sum_{\substack{c_{ir}\in C_{rft}\\ r\in \sR}}\left({t + T^\Delta}\right) - \bigM \sum_{\substack{c_{ir}\in C_{rft}\\ r\in \sR}}{\left( 2-{q_{c_{ir}fk}}-x_{{c_{ir}fkt}}^{\alpha}\right)} \leq \sum_{\substack{c_{ir}\in C_{rft}\\ r\in \sR}}{d_{c_{ir}}} + \sum_{\substack{c_{ir}\in C_{rft}\\ r\in \sR}}{T_{c_{ir}f}^\tau{q_{c_{ir}fk}}} \\+\sum_{\substack{c_{ir}\in C_{rft}\\ r\in \sR}}{{w_{c_{ir}fk}}},
  \end{align*}
  where $\sum_{\substack{c_{ir}\in C_{rft}\\ r\in \sR}}{x_{c_{ir}fkt}^{\alpha}} = 0$. 
  If we assume $T^\Delta$ is small ($T^\Delta \to 0$), then we can eliminate
  $T^\Delta$ and $\epsilon$. 
  Based on $\sum_{\substack{c_{ir}\in C_{rft}\\ r\in \sR}}{d_{c_{ir}}} + \sum_{\substack{c_{ir}\in C_{rft}\\ r\in \sR}}{w_{c_{ir}fk}}+ \sum_{\substack{c_{ir}\in C_{rft}\\ r\in \sR}}{T_{c_{ir}f}^\tau {q_{c_{ir}fk}}} \leq \sum_{\substack{c_{ir}\in C_{rft}\\ r\in \sR}}{t}$ and $\sum_{\substack{c_{ir}\in C_{rft}\\ r\in \sR}}{d_{c_{ir}}} + \sum_{\substack{c_{ir}\in C_{rft}\\ r\in \sR}}{{w_{c_{ir}fk}}} + \sum_{\substack{c_{ir}\in C_{rft}\\ r\in \sR}}{T_{c_{ir}f}^\tau{q_{c_{ir}fk}}}$ $\geq \sum_{\substack{c_{ir}\in C_{rft}\\ r\in \sR}}{t} - \bigM$ we can conclude $\sum_{\substack{c_{ir}\in C_{rft}\\ r\in \sR}}{w_{c_{ir}fk}}$ $\leq \sum_{\substack{c_{ir}\in C_{rft}\\ r\in \sR}}{t} - \sum_{\substack{c_{ir}\in C_{rft}\\ r\in \sR}}{d_{c_{ir}}} - \sum_{\substack{c_{ir}\in C_{rft}\\ r\in \sR}}{T_{c_{ir}f}^\tau {q_{c_{ir}fk}}}$. Since $t$ is a non-negative variable and the model minimizes the $w_{c_{ir}fk}$,  $w_{c_{ir}fk} = 0$. 
\end{proof}

\clearpage
\bibliographystyle{elsarticle-harv} 
\bibliography{our_bib}
\newpage 

\renewcommand{\thesection}{Supplementary Materials}
\section{Genetic algorithm pseudocode details}\label[app]{gaappendix}
\renewcommand{\thetable}{S.\arabic{table}}
\renewcommand{\thefigure}{S.\arabic{figure}}
\renewcommand{\thealgocf}{S.\arabic{algocf}}
\setcounter{table}{0}
\setcounter{figure}{0}
\setcounter{algocf}{0}

\begin{algorithm}[!htb]
\SetArgSty{textnormal}
\scriptsize
\SetKwInOut{Input}{Input}
\SetKwInOut{Output}{Output}
\Input{~ {$\sC_r$, $\sF$, $\sF_{c_{ir}}$, $\sK$, $\sR$, $\overline{B}$, $B^\iota$, $B^\omega$, $B_r^\iota$, $B_r^\omega$, $C^{\alpha=0}$, $C^\rho$, $C^\xi_k$, $C_k^\nu$, $C_f^\phi$, $\barbelow{N}$, $\bar{N}$, $N^{pop}$, $N^{iter}$, $N^{parents}$, $P^{mutate}$, $R_k$, $\overline{T}$, $T^\Delta$, $T^\rho_r$, $T^\mu_r$, $T_{c_{ir}c_{jr}}^\tau$, $T_{c_{ir}f}^\delta$, $T_{c_{ir}}^\kappa$}}
\Output{~ {$Sols(1)$}}
\SetKwFunction{FGAMain}{\textproc{GAMain}}
  \SetKwProg{Fn}{Function}{:}{\KwRet {$Sols(1)$}}
  \Fn{\FGAMain{}}
  {$Sols\gets \{\}$\;
\For{$i\gets 1$ \KwTo $N^{pop}$,}{
$Sols(i) \gets$\textproc{Initialization}($\sC_r$, $\sF_{c_{ir}}$, $\sK$, $\sR$, $\overline{B}$, $B^\iota$, $B^\omega$, $C^{\alpha=0}$, $T^\mu_r$, $T_{c_{ir}f}^\delta$)\;
\tcc{Create $N^{pop}$ solutions using \textproc{Initialization} function.}
}
\For{$iter \gets  1$ \KwTo $N^{iter}$,}{
\For{$i \gets  1$ \KwTo $N^{pop}$,}{
$\textbf{C}(Sols(i))\gets$ \textproc{Evaluator}($\sC_r$, $\sF$, $\sK$, $\sR$, $\sT$, $\overline{B}$, $B_r^\iota$, $B_r^\omega$, $C^\rho$, $C^\xi_k$, $C_k^\nu$, $C_f^\phi$, $\barbelow{N}$, $\bar{N}$, $R_k$, $T_{c_{ir}c_{jr}}^\tau$, $T_{c_{ir}}^\kappa$, $T_{c_{ir}f}^\delta$, $T_r^\rho$, $Sols(i)$)\;}
Sort $Sols$ ascending based on $\textbf{C}(Sols)$ and set $Sols$ to the first $N^{parents}$ of $Sols$\;
$CrossoverSols\gets \{\}$\;
\For{$i\gets 1$ \KwTo $N^{parents}$,}{\For{$j \gets  1$ \KwTo $N^{parents}$,}{\If{$i\neq j$,}{
$CrossoverSols(i)\gets$ \textproc{Crossover}($Sols$, $i$, $j$)\;
}}}
$MutationSols\gets \{\}$\;
{\For{$i \gets  1$ \KwTo $N^{parents}$,}{$MutationSols(i)\gets$ \textproc{Mutation}($Sols(i)$, $P^{mutate}$)\;
}}
\For{$i \gets  1$ \KwTo $N^{parents}$,}
{$\textbf{C}(CrossoverSols(i))\gets$ \textproc{Evaluator}($\sC_r$, $\sF$, $\sK$, $\sR$, $\sT$, $\overline{B}$, $B_r^\iota$, $B_r^\omega$, $C^\rho$, $C^\xi_k$, $C_k^\nu$, $C_f^\phi$, $\barbelow{N}$, $\bar{N}$, $R_k$, $T_{c_{ir}c_{jr}}^\tau$, $T_{c_{ir}}^\kappa$, $T_{c_{ir}f}^\delta$, $T_r^\rho$, $CrossoverSols(i)$)\;
$\textbf{C}(MutationSols(i))\gets$ \textproc{Evaluator}($\sC_r$, $\sF$, $\sK$, $\sR$, $\sT$, $\overline{B}$, $B_r^\iota$, $B_r^\omega$, $C^\rho$, $C^\xi_k$, $C_k^\nu$, $C_f^\phi$, $\barbelow{N}$, $\bar{N}$, $R_k$, $T_{c_{ir}c_{jr}}^\tau$, $T_{c_{ir}}^\kappa$, $T_{c_{ir}f}^\delta$, $T_r^\rho$, $MutationSols(i)$)\;} 
Merge $Sols$, $CrossoverSols$, and $MutationSols$ to form $newSols$\;
Sort $newSols$ ascending based on $\textbf{C}(newSols)$\;
Replace $Sols$ with first $N^{pop}$ of $newSols$;}}
 \caption{Pseudocode for genetic algorithm}\label[algo]{GA_main}
\end{algorithm}

\begin{algorithm}[H]
\SetArgSty{textnormal}
\SetAlgoLined
\scriptsize
\DontPrintSemicolon
  \SetKwFunction{FAbs}{\textproc{Abs}}
  \SetKwProg{Fn}{Function}{:}{\KwRet {absolute value of $a$}}
  \Fn{\FAbs{$a$}}{
        \tcc{finds absolute value of $a$.}
        }
  \;
   \SetKwFunction{FEnumerate}{\textproc{Enumerate}}
  \SetKwProg{Fn}{Function}{:}{\KwRet {index and item of $\sA$}}
  \Fn{\FEnumerate{$\sA$}}{
        \tcc{returns index and item of $\sA$.}
        }
  \;
  
  \SetKwFunction{FInt}{\textproc{Int}}
  \SetKwProg{Fn}{Function}{:}{\KwRet {integer part of $a$}}
  \Fn{\FInt{$a$}}{
        \tcc{finds integer part of $a$.}
        }
  \;
    
  \SetKwFunction{FRand}{\textproc{Rand}}
  \SetKwProg{Fn}{Function}{:}{\KwRet {$c$, such that $a\leq c\in\mathbb{R} \leq b$}}
  \Fn{\FRand{$a$, $b$}}{
        \tcc{finds a pseudo-random real number $c$ between $a$ and $b$.}
        }\;
  \SetKwFunction{FRandChoice}{\textproc{RandChoice}}
  \SetKwProg{Fn}{Function}{:}{\KwRet {$\sZ$, that is $\sZ\subseteq\sA$ and $|\sZ|=n$}}
  \Fn{\FRandChoice{$\sA$, $\sB$ (optional), $n$}}{
        \tcc{choose $n$ elements from list (or set) $\sA$ given a list of weights $\sB$.}
        }
  \;
  \SetKwFunction{FSortGiven}{\textproc{SortGiven}}
  \SetKwProg{Fn}{Function}{:}{\KwRet {$\sA$ sorted}}
  \Fn{\FSortGiven{$\sA$, $\sB$, $R$ (optional)}}{
        \tcc{sorts $\sA$ descending given $\sB$ if $R = True$, otherwise sorts $\sA$ ascending given $\sB$.}
        }
 \caption{Preliminary functions (see e.g., \cite{numpy} for details)}\label[algo]{Prelim_fn}
\end{algorithm}

\begin{algorithm}[H]
\SetArgSty{textnormal}
\scriptsize
\SetKwInOut{Input}{Input}
\SetKwInOut{Output}{Output}
\Input{~{$\sC_r$, $\sF_{c_{ir}}$, $\sK$, $\sR$, $\overline{B}$, $B^\iota$, $B^\omega$, $C^{\alpha=0}$, $T^\mu_r$, $T_{c_{ir}f}^\delta$}}
\Output{~{$\widehat{\sF}_{c_{ir}}$, $\widehat{\sK}_{c_{ir}}$}}
  \SetKwFunction{FInitialization}{\textproc{Initialization}}
  \SetKwProg{Fn}{Function}{:}{\KwRet {$\widehat{\sF}_{c_{ir}}$, $\widehat{\sK}_{c_{ir}}$}}
  \Fn{\FInitialization{}}
  {\tcc{Parameter $C^{\alpha=0}$ is the probability that a recharging will not occur after servicing a given cluster, and $T^\rho_r$ is the travel time for route $r$ excluding service and recharging times.}
  $f^{closest}\gets \{\}$; $\widehat{\sF}_{c_{ir}}\gets \{\}$\; \tcc{$\widehat{\sF}_{c_{ir}}$ 
  stores the facility that will be visited after serving $c_{ir}$; visiting $0$ denotes no recharging.}
  \For{$r\in \sR$,}{\For{$c_{ir}\in\sC_r$,}{
        $W\gets$;\ \tcp{a list of weights in $[0,1]$ based on $\sF_{c_{ir}}$}
        $f^{closest}(c_{ir})\gets$\textproc{RandChoice}(\textproc{SortGiven}($\sF_{c_{ir}}, T_{c_{ir}f}^\delta)$, $W$, 1)\;
        $\widehat{\sF}_{c_{ir}}\gets$\textproc{RandChoice}($[0,f^{closest}(c_{ir})]$, $[C^{\alpha=0},~1-C^{\alpha=0}]$, 1)\;
        }}
  \For{$r\in \sR$,}{\If{$\widehat{\sF}_{c_{ir}}=0~\forall c_{ir}\in\sC_r$,}
  {$\sC^{select}\gets$\textproc{RandChoice}($\sC_r$, $\ceil*{(\textproc{Abs}((B^\iota_r-B^\omega_r-T^\mu_r)/\overline{B})))}$)\;
  \tcc{$\ceil{(\textproc{Abs}((B^\iota_r-B^\omega_r-T^\mu_r)/\overline{B})))}$ denotes the minimum number of charging facilities need to be visited.}}
  \For{$c^\prime\in\sC^{select}$,}{$\widehat{\sF}_{c_{ir}}(c^\prime )\gets f^{closest}(c^\prime)$\;}}
  $\widehat{\sK}_{c_{ir}}\gets \{\}$\;\tcc{$\widehat{\sK}_{c_{ir}}$ 
  stores the charger type that will be used after $c_{ir}$; $0$ refers to not visiting a recharging station.}
  $W\gets$;\ \tcp{a list of weights in $[0,1]$ based on charger type}
  \For{$r\in\sR$,}{\For{$c_{ir}\in\sC_r$,}{
  \uIf{$\widehat{\sF}_{c_{ir}} \neq 0$,}{
  $\widehat{\sK}_{c_{ir}} \gets$\textproc{RandChoice}($\sK$, $W$, 1)\;
  }
  \Else{$\widehat{\sK}_{c_{ir}} \gets 0$\;}
  }}
  }
 \caption{Pseudocode for \textproc{Initialization}}\label[algo]{Init_pseudo}
\end{algorithm}

\begin{algorithm}[H]
\SetArgSty{textnormal}
\scriptsize
\SetKwInOut{Input}{Input}
\SetKwInOut{Output}{Output}
\Input{~ {$Sols$, $i$, $j$}}
\Output{~ {$newSols$}}
\SetKwFunction{FCrossover}{\textproc{Crossover}}
  \SetKwProg{Fn}{Function}{:}{\KwRet {$newSols$}}
  \Fn{\FCrossover{}}
  {$newSols\gets \{\}$\;
  \For{$element\gets 1$ \KwTo $|Sols(i)|$,}
  {\uIf{\textproc{RandChoice}($[i,j]$, $1$) = $i$,}
  {$newSols(element)\gets Sols(i)(element)$\;}
  \Else{$newSols(element)\gets Sols(j)(element)$\;}}
  }
 \caption{Pseudocode for \textproc{Crossover}}\label[algo]{crossover_pseudo}
\end{algorithm}

\begin{algorithm}[H]
\SetArgSty{textnormal}
\scriptsize
\SetKwInOut{Input}{Input}
\SetKwInOut{Output}{Output}
\Input{~ {$Sols$, $P^{mutate}$}}
\Output{~ {$newSols$}}
\SetKwFunction{FMutation}{\textproc{Mutation}}
  \SetKwProg{Fn}{Function}{:}{\KwRet {$newSols$}}
  \Fn{\FMutation{}}
  {$N^{mutate} \gets \textproc{Int}(|Sols| P^{mutate})$\;
  \tcc{The percent of solution to be mutated is denoted by $P^{mutate}$, and $N^{mutate}$ is the number of elements to mutate.}
  $newSols\gets \{\}$\;
  \For{$i\in$\textproc{RandChoice}($Sols$, $N^{mutate}$),}
  {\uIf{$Sols(i)$ indicates not recharging after visiting $i$,}
  {$newSols(i)\gets$; \tcp{Closest facility in terms of travel time and charger type $1$}}
  \Else{$newSols(i)\gets$; \tcp{Next closest facility in terms of travel time and charger type incremented by $1$}}}
  }
 \caption{Pseudocode for \textproc{Mutation}}\label[algo]{mutation_pseudo}
\end{algorithm}

\begin{algorithm}[H]
\SetArgSty{textnormal}
\scriptsize
\SetKwInOut{Input}{Input}
\SetKwInOut{Output}{Output}
\Input{~ {$\sC_r$, $\sF$, $\sK$, $\sR$, $\sT$, $\overline{B}$, $B_r^\iota$, $B_r^\omega$, $C^\rho$, $C^\xi_k$, $C_k^\nu$, $C_f^\phi$, $\barbelow{N}$, $\bar{N}$, $R_k$, $T_{c_{ir}c_{jr}}^\tau$, $T_{c_{ir}}^\kappa$, $T_{c_{ir}f}^\delta$, $T_r^\rho$, $q_{c_{ir}fk}$}}
\Output{~ {$\textbf{C}$}}
\SetKwFunction{FEvaluator}{\textproc{Evaluator}}
  \SetKwProg{Fn}{Function}{:}{\KwRet {$\textbf{C}$}}
  \Fn{\FEvaluator{}}
{
$\bar{z}_{fk} \gets \sum_{{c_{ir}}\in\sC_{rft}}{q_{c_{ir}fk}}$\;
\For{$f \in \sF$,} {
\For{$k \in \sK$,} {
\If{$\sum_{{c_{ir}}\in\sC_{rft}}{q_{c_{ir}fk}} \geq 0$,}{
$\barbelow{z}_{fk}\gets 1$\;}
}}
$(\barbelow{b}_{c_{ir}}, \barbelow{b}_{c_{ir}f}^\prime, \barbelow{\textbf{C}}, \barbelow{d}_{c_{ir}}, \barbelow{u}_{c_{ir}fk}, 
\barbelow{w}_{c_{ir}fk}, \barbelow{x}_{c_{ir}fkt}, \barbelow{y}_f) \gets$  \textproc{LowerLevelEvaluator}($\sC_r$, $\sF$, $\sK$, $\sR$, $\sT$, $\overline{B}$, $B_r^\iota$, $B_r^\omega$, $C^\rho$, $C^\xi_k$, $C_k^\nu$, $C_f^\phi$, $R_k$, $T_r^\rho$, $T_{c_{ir}c_{jr}}^\tau$, $T_{c_{ir}}^\kappa$, $T_{c_{ir}f}^\delta$, $q_{c_{ir}fk}$, $\barbelow{z}_{fk}$)\;
$(\bar{b}_{c_{ir}}, \bar{b}_{c_{ir}f}^\prime, \bar{\textbf{C}}, \bar{d}_{c_{ir}}, \bar{u}_{c_{ir}fk}, 
\bar{w}_{c_{ir}fk}, \bar{x}_{c_{ir}fkt}, \bar{y}_f) \gets$  \textproc{LowerLevelEvaluator}($\sC_r$, $\sF$, $\sK$, $\sR$, $\sT$, $\overline{B}$, $B_r^\iota$, $B_r^\omega$, $C^\rho$, $C^\xi_k$, $C_k^\nu$, $C_f^\phi$, $R_k$, $T_r^\rho$, $T_{c_{ir}c_{jr}}^\tau$, $T_{c_{ir}}^\kappa$, $T_{c_{ir}f}^\delta$, $q_{c_{ir}fk}$, $\Bar{z}_{fk}$)\;
$\bar{z}_{fk}\gets \sum_{{c_{ir}}\in\sC_{rft}, t\in\sT}{\bar{x}_{c_{ir}fkt}}$\; 
\While{$\barbelow{\textbf{C}}$ is infeasible,}
{$\sL\gets \{\}$;

\For{$f \in \sF$,} {
\For{$k \in \sK$,} {
\If{$\barbelow{z}_{fk} \geq 0$,}{
$\sL \gets  \sL \cup \{(f,k)\}$\;
}
}}
{\textproc{Sort}($\sL$, $\sum_{{c_{ir}}\in\sC_{rft}}{\barbelow{w}_{c_{ir}fk}}$, $R = True$)\; $N \gets$ $\textproc{Rand}\left(\barbelow{N}, \bar{N}\right)$\;
\For{$\left(f, k\right) \in \sL$,} {$\barbelow{z}_{fk}\gets \barbelow{z}_{fk}+1$;}
$\bar{N} \gets  1/2\left(\bar{N}+\barbelow{N}\right);$
}}
\uIf{$\Bar{\textbf{C}} \leq \barbelow{\textbf{C}}$,}{
$z_{fk} \gets$ \textproc{zuUpdater}($\sC_r$, $\sF$, $\sK$, $\sR$, $\sT$, $\overline{B}$, $B_r^\iota$, $B_r^\omega$, $C^\rho$, $C^\xi_k$, $C_k^\nu$, $C_f^\phi$, $\barbelow{N}$, $\bar{N}$, $R_k$, $T_r^\rho$, $T_{c_{ir}c_{jr}}^\tau$, $T_{c_{ir}}^\kappa$, $T_{c_{ir}f}^\delta$, $q_{c_{ir}fk}$, $\barbelow{z}_{fk}$, $\Bar{z}_{fk}$)\;
}
\Else{
$z_{fk}\gets$\textproc{zlUpdater}($\sC_r$, $\sF$, $\sK$, $\sR$, $\sT$, $\overline{B}$, $B_r^\iota$, $B_r^\omega$, $C^\rho$, $C^\xi_k$, $C_k^\nu$, $C_f^\phi$, $\barbelow{N}$, $\bar{N}$, $R_k$, $T_r^\rho$, $T_{c_{ir}c_{jr}}^\tau$, $T_{c_{ir}}^\kappa$, $T_{c_{ir}f}^\delta$, $q_{c_{ir}fk}$, $\barbelow{z}_{fk}$, $\Bar{z}_{fk}$)\;
}
$(b_{c_{ir}}, b_{c_{ir}f}^\prime, \textbf{C}, u_{c_{ir}fk}, d_{c_{ir}}, w_{c_{ir}fk}, x_{c_{ir}fkt}, y_f) \gets$  \textproc{LowerLevelEvaluator}($\sC_r$, $\sF$, $\sK$, $\sR$, $\sT$, $\overline{B}$, $B_r^\iota$, $B_r^\omega$, $C^\rho$, $C^\xi_k$, $C_k^\nu$, $C_f^\phi$, $R_k$, $T_r^\rho$, $T_{c_{ir}c_{jr}}^\tau$, $T_{c_{ir}}^\kappa$, $T_{c_{ir}f}^\delta$, $q_{c_{ir}fk}$, $z_{fk}$)\;
}
\caption{Pseudocode for \textproc{Evaluator}}\label[algo]{Eval_pseudo}
\end{algorithm}

\begin{algorithm}[H]
\SetArgSty{textnormal}
\scriptsize
\SetKwInOut{Input}{Input}
\SetKwInOut{Output}{Output}
\Input{~ {$\sC_r$, $\sF$, $\sK$, $\sR$, $\sT$, $\overline{B}$, $B_r^\iota$, $B_r^\omega$, $C^\rho$, $C^\xi_k$, $C_k^\nu$, $C_f^\phi$, $R_k$, $T_r^\rho$, $T_{c_{ir}c_{jr}}^\tau$, $T_{c_{ir}}^\kappa$, $T_{c_{ir}f}^\delta$, $q_{c_{ir}fk}$, $z_{fk}$}}
\Output{~ {$Out$}}
\SetKwFunction{FLowerLevelEvaluator}{\textproc{LowerLevelEvaluator}}
\SetKwProg{Fn}{Function}{:}{\KwRet {$Out$}}
\Fn{\FLowerLevelEvaluator{}}{
    $y_f \gets$ \textproc{yCalculation}$\left(\sF, ~\sK, ~z_{fk}\right)$\;
    $\left(\sL_{fk}, b_{c_{ir}f}^\prime, u_{c_{ir}fk}\right)\gets$ \textproc{uCalculation}$(\sC_r$, $\sF$, $\sK$, $\sR$, $\overline{B}$, $B_r^\iota$, $B_r^\omega$, $R_k$, $T_r^\rho$, $T_{c_{ir}c_{jr}}^\tau$, $T_{c_{ir}}^\kappa$, $q_{c_{ir}fk}$\;
    $\left(b_{c_{ir}}, d_{c_{ir}}, d_{c_{ir}f}^\prime d_{c_{ir}f}^{\prime\prime}, w_{c_{ir}fk}\right) \gets$ \textproc{wCalculation}($\sC_r$, $\sF$, $\sK$, $\sL_{fk}$, $\sR$, $T_{c_{ir}f}^\delta$, $T_{c_{ir}c_{jr}}^\tau$, $T_{c_{ir}}^\kappa$, $q_{c_{ir}fk}$, $u_{c_{ir}fk}$)\;
    $x_{c_{ir}fkt} \gets$\textproc{xCalculation}$\left(\sC_r,~ \sF,~ \sK, ~\sT, ~d_{c_{ir}f}^\prime, ~d_{c_{ir}f}^{\prime\prime}, ~q_{c_{ir}fk}, ~u_{c_{ir}fk}\right)$\;
    $\textbf{C} \gets \sum_{\substack{c_{ir}\in \sC_r,r \in \sR,\\f\in \sF_{c_{ir}}, k\in \sK}} \left[C^\rho\left(T_{c_{ir}f}^\delta q_{c_{ir}fk}+w_{c_{ir}fk}\right)+(C^\rho +C^\xi_k) u_{c_{ir}fk}\right] + \sum_{f\in\sF} C_f^\phi y_f + \sum_{f\in\sF,k\in\sK}C_k^\nu z_{fk}$\;
    \uIf{\labelcref{x_formulation}--\labelcref{operation_time_constraint} hold,}{$Out \gets$ ($b_{c_{ir}}$, $b_{c_{ir}f}^\prime$, $\textbf{C}$, $d_{c_{ir}}$, $u_{c_{ir}fk}$, $w_{c_{ir}fk}$, $x_{c_{ir}fkt}$, $y_f$);}\Else{$Out\gets None$;}
}
\caption{Pseudocode for \textproc{LowerLevelEvaluator}}\label[algo]{LLE_pseudo}
\end{algorithm}

\begin{algorithm}[H]
\SetArgSty{textnormal}
\scriptsize
\SetKwInOut{Input}{Input}
\SetKwInOut{Output}{Output}
\Input{~ {$\sC_r$, $\sF$, $\sK$, $\sR$, $\sT$, $\overline{B}$, $B_r^\iota$, $B_r^\omega$, $C^\rho$, $C^\xi_k$, $C_k^\nu$, $C_f^\phi$, $\barbelow{N}$, $\bar{N}$, $R_k$, $T_r^\rho$, $T_{c_{ir}c_{jr}}^\tau$, $T_{c_{ir}}^\kappa$, $T_{c_{ir}f}^\delta$, $q_{c_{ir}fk}$, $\barbelow{z}_{fk}$, $\Bar{z}_{fk}$}}
\Output{~ {$z_{fk}$}}
  \SetKwFunction{FzuUpdater}{\textproc{zuUpdater}}
  \SetKwProg{Fn}{Function}{:}{\KwRet {$z_{fk}$}}
  \Fn{\FzuUpdater{}}{
  \While{True,}{$\sL \gets \{\}$\;
  \For{$f \in \sF$,} {
  \For{$k \in \sK$,} {
  \If{$\bar{z}_{fk}-\barbelow{z}_{fk} \geq 1$ \text{and} $\barbelow{z}_{fk} \geq 0$ ,}{
  $\sL \gets  \sL \cup \{(f,k)\}$\;
  }
  }}
  \textproc{Sort}($\sL$, $\Bar{z}_{fk}$, $R = True$)\;
  $N \gets$ \textproc{Rand}$\left(\barbelow{N}, \bar{N}\right)$\;
  $\bar{\bar{z}}_{fk}\gets \Bar{z}_{fk}$\;
  \For{$\left(f, k\right) \in \sL$,} {$\bar{\bar{z}}_{fk}\gets \bar{z}_{fk}-1$;}
  $\bar{N} \gets  1/2\left(\bar{N}+\barbelow{N}\right)$\;
  $(\Bar{\bar{b}}_{c_{ir}}, \Bar{\bar{b}}_{c_{ir}f}^\prime, \Bar{\bar{\textbf{C}}}, \Bar{\bar{d}}_{c_{ir}}, \Bar{\bar{u}}_{c_{ir}fk}, 
  \Bar{\bar{w}}_{c_{ir}fk}, \Bar{\bar{x}}_{c_{ir}fkt}, \Bar{\bar{y}}_f) \gets$ \textproc{LowerLevelEvaluator}($\sC_r$, $\sF$, $\sK$, $\sR$, $\sT$, $\overline{B}$, $B_r^\iota$, $B_r^\omega$, $C^\rho$, $C^\xi_k$, $C_k^\nu$, $C_f^\phi$, $R_k$, $T_r^\rho$, $T_{c_{ir}c_{jr}}^\tau$, $T_{c_{ir}}^\kappa$, $T_{c_{ir}f}^\delta$, $q_{c_{ir}fk}$, $\Bar{\Bar{z}}_{fk}$)\;
  \uIf{$\bar{\bar{\textbf{C}}} \leq \bar{\textbf{C}}$,}{
  $\bar{z}_{fk}\gets \bar{\bar{z}}_{fk}$\;
  }\Else{$z_{fk}\gets \bar{z}_{fk}$\;
  $Break$\;}
  }
  }
  \caption{Pseudocode for calculating $\bar{z}_{fk}$}\label[algo]{z_upper_pseudo}
\end{algorithm}

\begin{algorithm}[H]
\SetArgSty{textnormal}
\scriptsize
\SetKwInOut{Input}{Input}
\SetKwInOut{Output}{Output}
\Input{~ {$\sC_r$, $\sF$, $\sK$, $\sR$, $\sT$, $\overline{B}$, $B_r^\iota$, $B_r^\omega$, $C^\rho$, $C^\xi_k$, $C_k^\nu$, $C_f^\phi$, $\barbelow{N}$, $\bar{N}$, $T_r^\rho$, $T_{c_{ir}c_{jr}}^\tau$, $T_{c_{ir}}^\kappa$, $T_{c_{ir}f}^\delta$, $R_k$, $q_{c_{ir}fk}$, $\barbelow{z}_{fk}$, $\Bar{z}_{fk}$}}
\Output{~ {$z_{fk}$}}
  \SetKwFunction{FzlUpdater}{\textproc{zlUpdater}}
  \SetKwProg{Fn}{Function}{:}{\KwRet {$z_{fk}$}}
  \Fn{\FzlUpdater{}}{
  \While{True,}{$\sL \gets  \{\}$\;
  \For{$f \in \sF$,} {
  \For{$k \in \sK$,} {
  \If{$\bar{z}_{fk}-\barbelow{z}_{fk} \geq 1~\text{and}~\barbelow{z}_{fk} \geq 0$,}{
  $\sL \gets  \sL \cup \{(f,k)\}$\;
  }}}
  \textproc{Sort}($\sL$, $\sum_{{c_{ir}}\in\sC_{rft}}{\barbelow{w}_{c_{ir}fk}}$, $R = True$)\;
  $N\gets$ $\textproc{Rand}\left(\barbelow{N}, \bar{N}\right)$\;
  $\barbelow{\barbelow{z}}_{fk}\gets \barbelow{z}_{fk}$\;
  \For{$\left(f, k\right) \in \sL$,} {$\barbelow{\barbelow{z}}_{fk}\gets \barbelow{\barbelow{z}}_{fk}+1$\;}
  $\bar{N} \gets  1/2\left(\bar{N}+\barbelow{N}\right)$\;
  $(\barbelow{\barbelow{b}}_{c_{ir}}, \barbelow{\barbelow{b}}_{c_{ir}f}^\prime, \barbelow{\barbelow{\textbf{C}}}, \barbelow{\barbelow{d}}_{c_{ir}}, \barbelow{\barbelow{u}}_{c_{ir}fk}, \barbelow{\barbelow{w}}_{c_{ir}fk}, \barbelow{\barbelow{x}}_{c_{ir}fkt}, \barbelow{\barbelow{y}}_f)
  \gets$ \textproc{LowerLevelEvaluator}($\sC_r$, $\sF$, $\sK$, $\sR$, $\sT$, $\overline{B}$, $B_r^\iota$, $B_r^\omega$, $C^\rho$, $C^\xi_k$, $C_k^\nu$, $C_f^\phi$, $R_k$, $T_r^\rho$, $T_{c_{ir}c_{jr}}^\tau$, $T_{c_{ir}}^\kappa$, $T_{c_{ir}f}^\delta$, $q_{c_{ir}fk}$, $\barbelow{\barbelow{z}}_{fk}$)\;
  \uIf{$\barbelow{\barbelow{\textbf{C}}} \leq \barbelow{\textbf{C}}$,}{
  $\barbelow{z}_{fk}\gets \barbelow{\barbelow{z}}_{fk}$\;
  }\Else{$z_{fk} \gets \barbelow{z}_{fk}$\;
  $Break$\;}
  }}
  \caption{Pseudocode for calculating $\barbelow{z}_{fk}$}\label[algo]{z_lower_pseudo}
\end{algorithm}

\begin{algorithm}[H]
\SetArgSty{textnormal}
\scriptsize
\SetKwInOut{Input}{Input}
\SetKwInOut{Output}{Output}
\Input{~ {$\sF$, $\sK$, $z_{fk}$}}
\Output{~ {$y_f$}}
\SetKwFunction{FyCalculation}{\textproc{yCalculation}}
\SetKwProg{Fn}{Function}{:}{\KwRet{$y_f$}}
\Fn{\FyCalculation{}}{
\For{$f \in \sF$,} {
\uIf{$\sum_{k\in\sK}{z_{fk}} \geq 0$,}{
$y_f \gets  1$\;
}
\Else{
$y_f \gets  0$\;
}
}
}

\caption{Pseudocode for calculating $y_f$}\label[algo]{y_pseudo}
\end{algorithm}

\begin{algorithm}[H]
\SetArgSty{textnormal}
\scriptsize
\SetKwInOut{Input}{Input}
\SetKwInOut{Output}{Output}
\Input{~ {$\sC_r$, $\sF$, $\sK$, $\sR$, $\overline{B}$, $B_r^\iota$, $B_r^\omega$, $T_r^\rho$, $T_{c_{ir}c_{jr}}^\tau$, $T_{c_{ir}}^\kappa$, $q_{c_{ir}fk}$}}
\Output{~ {$\sL_{fk}, b_{c_{ir}f}^\prime, u_{c_{ir}fk}$}}
\SetKwFunction{FuCalculation}{\textproc{uCalculation}}
\SetKwProg{Fn}{Function}{:}{\KwRet {$\sL_{fk}, b_{c_{ir}f}^\prime, u_{c_{ir}fk}$}}
\Fn{\FuCalculation{}}{
    \For{$r \in \sR$,} {
        \For{$f \in \sF$,}{\For{$k \in \sK$,}{$\sL_{fk} \gets  \{\}$\;}}
        $e_r \gets \textproc{Int}\left(\textproc{Abs}\left(B_r^\iota - B_r^\omega - T_r^\rho\right)+1\right)$\;
        $\sL \gets  \{\}$\;
        \For{$c_{ir}\in\sC_r$, }{
          \For{$f\in\sF$,}{
            \For{$k \in \sK$, }{
                \If{$q_{c_{ir}fk} = 1$}{$\sL \cup \{\left(c_{ir}fk\right)\}$\;
                $\sL_{fk} \cup \{\left(rc_{ir}\right)\}$\;}
            }
          }
        } 
        $m \gets  1$\;
        \While{$e_r \geq 0 $, }{
            \uIf{$m = 1$, }{
                $\left(c_{ir}fk\right)\gets \sL\{m\}$\;
                $b_{c_{ir}f}^\prime \gets  B_r^\iota - \sum_{c_{jr}\gets c_{0r}}^{c_{jr} \gets  c_{i-1,r}}{T_{c_{jr}c_{j+1,r}}^\tau} - T_{c_{ir}f}^\tau$\;
                $u_{c_{ir}fk} \gets  \min\{e_r, \overline{B} - b_{c_{ir}f}^\prime\}$\;
            }
            \Else{
                $\left(c_{ir}fk\right)\gets \sL\{m\}$\;
                $\left(c_{ir}^\prime f^\prime k^\prime\right)\gets \sL\{m-1\}$\;
                $b_{c_{ir}f}^\prime \gets  b_{c_{ir}^\prime f^\prime}^\prime + u_{c_{ir}^\prime f^\prime k^\prime} - \sum_{c_{jr}\gets c_{0r}}^{c_{jr} \gets  c_{i-1,r}}{T_{c_{jr}c_{j+1,r}}^\tau} - T_{c_{ir}f}^\tau$\;
                $u_{c_{ir}fk} \gets  \min\{e_r, \overline{B} - b_{c_{ir}f}^\prime\}$\;
            }
            $m \gets m+1$\;
            $e_r \gets e_r-u_{c_{ir}fk}$\;
        }
    }
}
\caption{Pseudocode for calculating $u_{c_{ir}fk}$}\label[algo]{u_pseudo}
\end{algorithm}

\begin{algorithm}[H]
\SetArgSty{textnormal}
\scriptsize
\SetKwInOut{Input}{Input}
\SetKwInOut{Output}{Output}
\Input{~ { $\sC_r$, $\sF$, $\sK$, $\sL_{fk}$, $\sR$, $T_{c_{ir}f}^\delta$, $R_k$, $T_{c_{ir}c_{jr}}^\tau$, $T_{c_{ir}}^\kappa$, $q_{c_{ir}fk}$, $u_{c_{ir}fk}$}}
\Output{~ {$b_{c_{ir}}, d_{c_{ir}}, d_{c_{ir}f}^\prime d_{c_{ir}f}^{\prime\prime}, w_{c_{ir}fk}$}}
\SetKwFunction{FwCalculation}{\textproc{wCalculation}}
\SetKwProg{Fn}{Function}{:}{\KwRet {$b_{c_{ir}}, d_{c_{ir}}, d_{c_{ir}f}^\prime d_{c_{ir}f}^{\prime\prime}, w_{c_{ir}fk}$}}
\Fn{\FwCalculation{}}{
    \For{$iter\gets 1$ \textbf{to} $|\sL_{fk}|$,}{
        \For{$r \in \sR$,} {
            $d_{c_{0r}} \gets 0$; $b_{c_{0r}} \gets B_r^\iota$\;
            \For{$c_{ir}\in\sC_r$, }{
              \uIf{$\sum_{f\in \sF, k\in\sK}{q_{c_{i-1,r}fk}} = 1$,}{$d_{c_{ir}} \gets d_{c_{i-1,r}} + T_{c_{i-1,r}f}^\delta + w_{c_{ir}fk} + u_{c_{ir}fk} + T_{c_{ir}}^\kappa$\;
              $b_{c_{ir}} \gets b_{c_{i-1,r}}-T_{c_{i-1,r}f}^\delta + R_k u_{c_{ir}fk}$ \;}
              \Else{$d_{c_{ir}} \gets d_{c_{i-1,r}} + T_{c_{i-1,r}c_{ir}}^\tau + T_{c_{ir}}^\kappa$ \;
              $b_{c_{ir}} \gets b_{c_{i-1,r}}-T_{c_{ir}c_{jr}}^\tau$ \;}
              \For{$f\in\sF$,}{
                \For{$k \in \sK$, }{
                    \If{$q_{c_{ir}fk} = 1$,}{
                        $d_{c_{ir}f}^\prime \gets d_{c_{ir}} + T_{c_{ir}f}^\tau$\;
                    }
                }
              }
            } 
        }
        \For{$f \in \sF$,}{\For{$k \in \sK$,}{\textproc{Sort}($\sL_{fk}$, $d_{c_{ir}f}^\prime$, $R = True$)\;
        \For{$j, \left(r, c_{ir}\right) \in \textproc{Enumerate}\left(\sL_{fk}\right)$,}{
            \uIf{$j = 0$,}{$w_{c_{ir}fk} \gets 0$\;
                $d_{c_{ir}f}^{\prime\prime}\gets d_{c_{ir}f}^{\prime}+u_{c_{ir}fk}$\;}
            \Else{$ \left(r^\prime, c_{ir^\prime}^\prime \right) \gets \sL_{fk}\{j\}$\;
                $w_{c_{ir}fk} \gets \max\{d_{c_{ir^\prime}^\prime f}^{\prime\prime} - d_{c_{ir}f}^{\prime},0\}$\;
                $d_{c_{ir}f}^{\prime\prime}\gets d_{c_{ir}f}^{\prime}+w_{c_{ir}fk}+u_{c_{ir}fk}$\;}
        }    
        }}
    }
}
\caption{Pseudocode for calculating $w_{c_{ir}fk}$}\label[algo]{w_pseudo}
\end{algorithm}

\begin{algorithm}[H]
\SetArgSty{textnormal}
\scriptsize
\SetKwInOut{Input}{Input}
\SetKwInOut{Output}{Output}
\Input{~ {$\sC_r$, $\sF$, $\sK$, $\sT$, $d_{c_{ir}f}^\prime$, $d_{c_{ir}f}^{\prime\prime}$, $q_{c_{ir}fk}$, $u_{c_{ir}fk}$}}
\Output{~ {$x_{c_{ir}fkt}$}}
\SetKwFunction{FxCalculation}{\textproc{xCalculation}}
\SetKwProg{Fn}{Function}{:}{\KwRet {$x_{c_{ir}fkt}$}}
\Fn{\FxCalculation{}}{
    \For{$c_{ir}\in\sC_r$,}{
        \For{$f\in \sF$,}{
            \For{$k\in\sK$,}{
                \For{$t\in\sT$,}{
                    \uIf{$t \geq d_{c_{ir}f}^\prime$ and $t \leq d_{c_{ir}f}^{\prime\prime}$,}{$x_{c_{ir}fkt} \gets 1$\;}
                    \Else{$x_{c_{ir}fkt} \gets 0$\;}
                }
            }
        }
    }
}
\caption{Pseudocode for calculating $x_{c_{ir}fkt}$}\label[algo]{x_pseudo}
\end{algorithm}

\vfill
\framebox{\parbox{.90\linewidth}{\scriptsize The submitted manuscript has been created by
        UChicago Argonne, LLC, Operator of Argonne National Laboratory (``Argonne'').
        Argonne, a U.S.\ Department of Energy Office of Science laboratory, is operated
        under Contract No.\ DE-AC02-06CH11357.  The U.S.\ Government retains for itself,
        and others acting on its behalf, a paid-up nonexclusive, irrevocable worldwide
        license in said article to reproduce, prepare derivative works, distribute
        copies to the public, and perform publicly and display publicly, by or on
        behalf of the Government.  The Department of Energy will provide public access
        to these results of federally sponsored research in accordance with the DOE
        Public Access Plan \url{http://energy.gov/downloads/doe-public-access-plan}.}}
\end{document}